\documentclass[fullpage,12pt]{article}
\small
\renewcommand{\baselinestretch}{1.2}
\usepackage{graphicx,epstopdf}
\usepackage{amsfonts}
\usepackage{amsmath}
\usepackage{mathrsfs}
\usepackage{amssymb}
\usepackage[T1]{fontenc}
\usepackage[latin1]{inputenc}
\usepackage{dsfont}
\usepackage{xcolor}
\usepackage{float}
\usepackage{enumitem}
\usepackage{bm}
\usepackage[justification=centering]{caption}

\usepackage[francais]{babel}

\newtheorem{theo}{Theorem}[section]
\newtheorem{rmk}{Remark}[section]
\newtheorem{lem}{Lemma}[section]
\newtheorem{lemma}{Lemma}
\newtheorem{prop}{Proposition}[section]
\newtheorem{corol}{Corollary}[section]

\begin{document}
\vskip 3mm
\noindent{\bf\large Asymptotic Normality of an M-estimator of regression function for truncated-censored data under $\alpha$-mixing condition}
\vskip 10mm
\noindent Hassiba BENSERADJ\footnote{Corresponding author}\\
\noindent Faculty of Sciences, University of Boumerdes UMBB\\
\noindent Avenue of independence, Boumerdes, Algeria\\
\noindent h.benseradj@univ-boumerdes.dz
\vskip 2mm
\noindent Zohra GUESSOUM\\
\noindent Lab. MSTD, Faculty of Mathematics, USTHB \\
\noindent BP 32, El-Alia 16111, Algeria\\
\noindent zguessoum@usthb.dz\\
\vskip 5mm
{{\bf \noindent ABSTRACT}}\\
\noindent In this paper, we establish weak consistency and asymptotic normality of an M-estimator of the regression function for left truncated and right censored (LTRC) model, where it is assumed that the observations form a stationary $\alpha$-mixing sequence. The result holds with unbounded objective function, and are applied to derive weak consistency and asymptotic normality of a kernel classical regression curve estimate. We also obtain a uniform weak convergence rate for the product-limit estimator of the lifetime and censored distribution under dependence, which are useful results for our study and other LTRC strong mixing framework. Some simulations are drawn to illustrate the results for finite sample.
\vskip 4mm
{{\bf \noindent MATHEMATICS SUBJECT CLASSIFICATION}}\\
{{62G05; 62G20}}
\vskip 8mm
\renewcommand{\baselinestretch}{1.2}\normalsize
\section{Introduction}
Let $\left( X_{i},Y_{i}\right) _{i=1:N}$ be a sequence of $N$ stationary random variable, identically distributed as $\left( X,Y\right) $, taking
value in the space $\mathbb{R}^{d}\times\mathbb{R}$. Our purpose is to study the interaction between $X$ and $Y$ which can be formulated as
\begin{equation*}
Y=m\left( X\right) +\epsilon ,
\end{equation*}%
where $m\left( .\right) $ is an unknown function, and $\epsilon $ is the random error variable independent of $X.$\\
The traditional approach to estimate $m\left( x\right) ,$ for a given $x\in\mathbb{R}^{d}$, applies the $L_{2}$ loss function, viz
\begin{equation*}
m\left( x\right) =\arg \underset{\theta }{\min }\mathbb{E}\left[ \left(Y-\theta \right) ^{2}|X=x\right] .
\end{equation*}
As is well known, the least squares method achieves optimum results when the residual distribution is gaussian. However, the method becomes unstable if
the noise has nonzero-mean, and /or if outliers are present in the data. To overcome this problem, one may replace the above square function by another one which assigns less weight on large value of $Y_{i},$ these methods are called robust regression, which leads to
\begin{equation*}
m\left( x\right) =\arg \underset{\theta }{\min }\mathbb{E}\left[ \rho \left(Y-\theta \right) |X=x\right] .
\end{equation*}
where $\rho \left( .\right) $ is a given real continuously differentiable loss function such as $\frac{d\rho }{d\theta }\left(.\right)=:\psi _{x}\left(.\right)$ is a strictly monotone function. Hence $m\left( x\right) $ exists, and can be seen as the unique solution with respect
to (w.r.t.) $\theta$ of
\begin{equation}\label{mx}
\mathbb{E}\left[ \psi _{x}\left( Y-\theta \right) |X=x\right] =0,
\end{equation}
and the corresponding non parametric kernel $M$-estimator $\widehat{m}\left(.\right) $ based on the sample size $N$ as a zero with respect to $\theta $
of
\begin{equation}\label{mNW}
\sum\limits_{i=1}^{N}K\left( \frac{X_{i}-x}{h_{N}}\right) \psi _{x}\left(Y_{i}-\theta \right) ,
\end{equation}
where $K\left( .\right) $ is a kernel function and $h_{N}$ is a smoothing parameter tending to zero along with $N$. Special cases of (\ref{mNW}) include the usual Nadaraya-Watson estimator. This motivation in combining the ideas of robustness with those of smoothing comeback to Tsybakov (1983) and Hardle (1984). Robinson (1984) studied this M-estimator in time series context and Collomb and Hardle (1986) established its uniform consistency with rates. Many other researches have been done for this type of M-estimator, see for example Bohente and Fraiman (1990) and the references therein for previous results and Laib and Ould Said (2000), Attaoui et al (2015) for recent advances and references.

However, all this researches consider the case when a complete data is available, while in most works where the survival time $Y$ is the variable
of interest, referred here as the lifetime, two different problems appear: The first one, if the time origin of the lifetime precedes the start of the
study, only subjects that fail after the beginning of the study are being followed, otherwise they are left truncated.  Wang and Liang (2012) defined a
new robust estimator adapted to this type of data based on the estimator of classical regression function proposed by Ould Said and Lemdani (2006) and established its weak and strong consistency (without rate), as well as its asymptotic normality for $\alpha $-mixing processes. The second
problem which can appear in survival data is right censoring, which arises when a subject leaves the study before an event occurs, or the study
ends before the event has occurred. For this type of data, Lemdani and Ould Said (2017) constructed a new estimator based on the so-called synthetic data and established its uniform strong consistency with rate and asymptotic normality using the Vapnik-Chervonenkis class (VC class) under independence.

Note that this two types of incomplete data may be occur simultaneously in a study, then the model is known as Left Truncated and Right Censored one
(LTRC). Recently, many results have been established for this type of data, we can cite Iglessias-P$\acute{e}$rez and Gonzalez-Manteiga (1999), Liang et al (2012) and liang and Liu (2013) for the conditional distribution and conditional probability density estimation. Liang et al (2015) and Guessoum and Tatachak (2020) for conditional quantile and hazard function estimation, respectively. In the context of nonparametric regression, Benseradj and Guessoum (2020) have constructed a kernel M-estimator and established its uniform strong consistency with rate under strong mixing condition.

The aim of this paper is to establish weak consistency and asymptotic normality for the $M$-estimator from truncated and censored data, under $%
\alpha $-mixing assumption, by considering weaker conditions that those given to get strong consistency. So we complete our first work on strong consistency, and we extend the result of Wang and Liang (2012) from truncated to censored and truncated data, by lightening the conditions on the mixing coefficient, and the result of Lemdani and Ould Said (2017) from the censored independent data to the censored dependent data (as particular case).

Our results hold with unbounded objective function $\psi_{x} $ and are applied to derive asymptotic normality of a kernel classical regression curve
estimate. We compare our results to those obtained by Liang (2011) in absence of censoring, and to those obtained by Guessoum and Ould Said (2012)
in absence of truncation. Recall that our resulting estimate of the regression curve $\left( \psi _{x}\left( u\right) =u\right) $ under censoring
is slightly different from that proposed by Carbonez et al (1995) (see Remark 4.3 in Benseradj and Guessoum (2020)). We show that our estimator has minimum asymptotical variance compared to the later.

The rest of the paper is organized as follows: In section 2, we give some definitions and notations related to our model. In section 3, we give the assumptions under which we state our main results. A simulation study that highlights our theoretical results is provided in section 4. Proofs and
auxiliary results are relegated to section 5. The Appendix section establishes a new supplementary results, concerning the asymptotic behavior of the Product-limit estimators for both lifetime and censoring time distributions, under $\alpha$-mixing.
\section{Model and Estimator}
Let $\left\{ \left( Y_{k},T_{k},W_{k}\right) ,1\leq k\leq N\right\} $ be a sequence of random vectors from $\left( Y,T,%
W\right) $, where $Y$ denotes the lifetime under study with continuous distribution function (d.f) $F.$ $T$ and $W$ are the variables of left
truncation and right censoring time with continuous (d.f's) $\ L$ and $G,$ respectively. Let%
\begin{equation*}
Z=\left( Y\wedge W\right) \text{ \ \ and \ \ }\delta =\mathbb{I}_{\left\{
Y\leq W\right\} ,}
\end{equation*}
where $\left( t\wedge u\right) :=\min (t,u),$ and $\delta $ is the indicator of censoring status. In random LTRC model one observe $\left( Z,T,\delta \right) $ only if $Z\geq T$. Set $\mu =\mathbb{P}\left( T\leq Z\right) ,$ then we need to assume that $\mu $ $>0,$ otherwise, nothing is
observable. Consider the presence of a covariate $X,$ and assume that $X$ admits a density $v.$ Then, denote by $\left( X_{i},Z_{i},T_{i},\delta _{i}\right),$ $i=1,2,...,n;$ $\left( n\leq N\right)$ a stationary random sample from $\left( X,Z,T,\delta \right) $ which one really
observe (ie, $T_{i}\leq Z_{i}$). In all the remaining of this paper we suppose that $T,W$ and $Y$ are mutually independent, then \ $Z$ has (d.f)
\begin{equation*}
H=1-\left( 1-F\right) \left( 1-G\right).
\end{equation*}
Note that $H$ can be estimated by the well known Lynden-Bell estimator $H_{n}$, since the random variable $Z$ verify the random left truncated (RLT) model
\begin{equation*}
H_{n}\left( y\right) {\LARGE =}1-\prod\limits_{i=1}^{n}\left( 1-\dfrac{\mathbb{I}_{\left\{ Z_{i}\leq y\right\} }}{nC_{n}\left(Z_{i}\right) }\right),
\end{equation*}
where $C_{n}\left( y\right)=\frac{1}{n}\sum\limits_{j=1}^{n}\mathbb{I}_{\left\{ T_{j}\leq y\leq Z_{j}\right\} }$, is the empirical estimator of $C\left( y\right) :=\mathbb{P}\left( T\leq y\leq Z|T\leq Z\right)$. Denote by $\mathbb{P}$ and $\mathbf{P}$ the probability measure related to the $N$ sample, and
the actually observed $n$ sample, respectively. Also $\mathbb{E}$ and $\mathbf{E}$ the expectation operators related to $\mathbb{P}$ and $\mathbf{P}
$ respectively. In the sequel, for any (d.f) $Q$, we denote its left and right support endpoints by $a_{Q}:=\inf \left\{ x:Q\left( x\right) >0\right\} $ and $b_{Q}:=\sup \left\{ x:Q\left( x\right) <1\right\} ,$ respectively. It is clear that $a_{H}=a_{F}\wedge a_{G},$ $b_{H}=b_{F}\wedge b_{G}.$ Tsai et al (1987) defined the product-limit (PL) estimator $F_{n}$ of the distribution function $F$ for the LTRC data, by
\begin{equation*}
1{\LARGE -}F_{n}\left( y\right) {\LARGE =}\prod\limits_{i=1}^{n}\left( 1-\dfrac{\mathbb{I}_{\left\{ Z_{i}\leq y\right\} }\delta_{i}}{nC_{n}\left(Z_{i}\right) }\right),
\end{equation*}
(which is called the TJW product-limit estimator). Note that, for strong mixing data, the uniform strong consistency of $F_{n}$
has been established, separately, in absence of truncation $(T=0)$ by Cai (1998) and in absence of censoring $(W=\infty )$ by Sun and Zhou (2001).
Chen and Day (2003) showed that the result of Sun and Zhou (2001) can be extended to both truncated and censored data. Uniform weak consistency of $F_{n}$ has been established, separately, for censoring data by Liang and Alvarez (2011) and for truncated data by Liang (2011). We extend their result in lemma \ref{LemmaA5} in Appendix for both truncated and censored data, which allowed us to deduce the uniform weak consistency with rate of $G_{n},$ the concomitant TJW estimator of the censored distribution $G,$ defined below in (\ref{G}), which are both interesting results for our study and other strong mixing framework.\\
Conditionally on the value of $n$ we define the joint distribution function of $Z$ and $T$ as
\begin{eqnarray*}
\mathbf{M}\left( z,t\right) :=\mathbf{P}\left( Z\leq z,T\leq t\right)
=\mu ^{-1}\int\limits_{-\infty}^{z}L\left( t\wedge u\right) dH\left( u\right) ,
\end{eqnarray*}
which implies that the conditional distribution functions of $Z$ and $T$ are defined respectively by
\begin{equation*}
\mathbf{H}\left( z\right) :=\mathbf{M}\left( z,\infty \right) =\mu
^{-1}\int\limits_{-\infty}^{z}L\left( u\right) dH\left( u\right), \text{ and }\text{\ } %
\mathbf{L}\left( t\right): =\mathbf{M}\left( \infty ,t\right) =\mu
^{-1}\int\limits_{-\infty}^{\infty }L\left( t\wedge u\right) dH\left( u\right) .
\end{equation*}
Of course, $\mathbf{H}\left( z\right) $ and $\mathbf{L}\left( t\right) $ can be empirically estimated by
\begin{equation*}
\mathbf{H}_{n}\left( z\right) :=\frac{1}{n}\sum\limits_{i=1}^{n}\mathbb{I}%
_{\left\{ Z_{i}\leq z\right\} },\text{ and }\text{ }\mathbf{L}_{n}\left( t\right): =\frac{%
1}{n}\sum\limits_{i=1}^{n}\mathbb{I}_{\left\{ T_{i}\leq t\right\}}.
\end{equation*}
Note that,
\begin{eqnarray*}
C\left( y\right) =\mathbb{P}\left( T\leq y\leq Z|T\leq Z\right)=\mu ^{-1}L\left( y\right) \left( 1-H\left( y\right) \right)
\end{eqnarray*}
since $H$ is continuous. Then
\begin{equation*}
\mu =\frac{L\left( y\right) \left( 1-F\left( y\right) \right) \left(1-G\left( y\right) \right) }{C\left( y\right) }, \text{  for } a_{H}\leq y<b_{H},
\end{equation*}
and following the idea of He and Yang (1998) a consistent estimate of $\mu $ is given by
\begin{equation}\label{2.2}
\mu _{n}=\frac{L_{n}\left( y\right) \left( 1-F_{n}\left( y\right) \right)
\left( 1-G_{n}\left( y\right) \right) }{C_{n}\left( y\right) },\text{ for
all }y\text{ such that }C_{n}\left( y\right) \neq 0,
\end{equation}
where $G_{n}$, the concomitant TJW estimator of the distribution function $G,$ is defined by\\
\begin{equation}\label{G}
1-G_{n}(y)=:\overline{G}_{n}\left(y\right)=\prod\limits_{i=1}^{n}\left( 1-\frac{\mathbb{I}_{\left\{Z_{i}\leq y\right\} }(1-\delta _{ i })}{nC_{n}\left( Z_{i}\right) }\right),
\end{equation}
and $L_{n}$ is the Lynden-Bell estimator of the distribution function $L$ defined by
\begin{equation}\label{L}
   L_{n}\left( y\right) =\prod\limits_{i:T_{i}>y}\left( 1-\dfrac{1}{nC_{n}\left( T_{i}\right) }\right).
\end{equation}
For RLT data, He and Yang (1998) proved that $\mu_{n}$ does not depend on $y$ and its value can be obtained for any $y$ such that $C_{n}(y)\neq 0$.
\vskip 4mm
\begin{rmk}\label{rem0}
Note that $\mu_{n}$ can be expressed as $ \mu _{n}=\frac{L_{n}\left( y\right)\left( 1-H_{n}\left( y\right)\right) }{C_{n}\left( y\right) }$, (since we have {\footnotesize $1-H_{n}\left( y\right)= \left( 1-F_{n}\left( y\right) \right)\left( 1-G_{n}\left( y\right) \right)$}), which is the same definition given by He and Yang (1998) with respect to the variable $Z$. Hence any asymptotic result for $\mu_{n}$ given for RLT data remains true under LTRC data.
 \end{rmk}
Under the current model, for LTRC data, Gijbels and Wang (1993) pointed out that $F$ can be estimated only if $a_{L}\leq a_{H},$ $%
b_{L}\leq b_{H}$ and following Ould Said and Lemdani (2006), these are necessary but not sufficient identifiability conditions. Throughout this paper, we assume that $a_{L}<a_{H},\text{ }b_{L}\leq b_{H}$, and
\begin{equation}\label{2.1}
\left( T\text{,}W\right) \text{ are independent of }\left( X,Y\right) .
\end{equation}
Define also the sub-conditional distribution function of $\left( Z,\delta =1,X\right)$ as
\begin{eqnarray*}
\mathbf{H}_{\mathbf{1}}\left( x,y\right) &:=&
\mathbf{P}\left( Z\leq y,\delta=1,X\leq x\right)  \notag \\
&=&\mathbb{P}\left( Z\leq y,\delta =1,X\leq x|T\leq Z\right)  \notag \\
&=&\mu ^{-1}\int\limits_{s\leq x}\int\limits_{a_{H}}^{y}L\left( t\right)\left( 1-G\left( t\right) \right) f_{X,Y}\left( s,t\right) dtds,
\end{eqnarray*}
where $f_{X,Y}\left( .,.\right) $ is the joint density function of $(X,Y)$.
We deduce that
\begin{equation}\label{2.4}
d\mathbf{H}_{1}\left( x,y\right) =\mu ^{-1}L\left( y\right) \left( 1-G\left(y\right) \right) f_{X,Y}\left( x,y\right) dxdy.
\end{equation}
In regression analysis, one expects to identify the relationship between $X$ and $Y$ via a robust regression. The robust regression function denoted $m\left(x\right) $, is implicitly defined as a solution of equation (\ref{mx}). The left hand side of (\ref{mx}) can be written as
\begin{equation*}
\mathbb{E}\left[ \psi _{x}\left( Y-\theta \right) |X=x\right] =\dfrac{\int\psi _{x}\left( y-\theta \right) f_{X,Y}\left( x,y\right) dy}{v\left(
x\right) },
\end{equation*}%
where $v\left( x\right) >0.$ Set%
\begin{equation*}
\Psi _{x}\left( x,\theta \right) :=\mathbb{E}\left[ \psi _{x}\left( Y-\theta
\right) |X=x\right] v\left( x\right) ,
\end{equation*}
then $m\left( x\right) $ can be viewed as a solution w.r.t $\theta$ of $\Psi _{x}\left(x,\theta \right) =0.$ Hence a natural estimator of $m\left( x\right) $ denoted $\widehat{m}\left(x\right) $ is defined as a zero w.r.t $\theta $ of
\begin{equation}\label{2.6}
\widehat{\Psi }_{x}\left( x,\theta \right) =0,
\end{equation}
where
\begin{equation}\label{2.5}
\widehat{\Psi }_{x}\left( x,\theta \right) :=\frac{\mu _{n}}{nh_{n}^{d}}%
\sum\limits_{i=1}^{n}K\left( \frac{x-X_{i}}{h_{n}}\right) \frac{\delta
_{i}\psi _{x}\left( Z_{i}-\theta \right) }{L_{n}\left( Z_{i}\right) \overline{G}_{n}\left( Z_{i}\right) },
\end{equation}
for all $Z_{i}$ such that $L_{n}\left( Z_{i}\right)\overline{G}_{n}\left( Z_{i}\right)\neq0$. The construction and the strong uniform consistency of the estimator defined in (\ref{2.5}) under strong mixing condition are detailed in Benseradj and Guessoum (2020). Another interesting intermediate (pseudo)-estimate of $\Psi _{x}$, which is a useful tool for establishing our results, is defined as
\begin{equation}\label{Psy}
\widetilde{\Psi }_{x}\left( x,\theta \right) :=\frac{\mu }{nh_{n}^{d}}%
\sum\limits_{i=1}^{n}K\left( \frac{x-X_{i}}{h_{n}}\right) \frac{\delta
_{i}\psi _{x}\left( Z_{i}-\theta \right) }{L\left( Z_{i}\right)\overline{G}\left( Z_{i}\right)}.
\end{equation}
for all $Z_{i}$ such that $a_{H}\leq Z_{i}<b_{H}$. In all the remaining of this paper, $\left\{ \left( X_{i},Z_{i},T_{i},\delta _{i}\right) ,\text{ }1\leq i\leq n\right\} $ is assumed to be a stationary $\alpha $-mixing sequences of random vectors. Recall that a sequence $\left\{ \zeta _{k},k\geq 1\right\} $ is said to be $\alpha $-mixing (strongly mixing) if the mixing coefficient
\begin{equation*}
\alpha \left( n\right) \overset{def}{:=}\sup\limits_{k\geq 1}\sup \left\{
\left\vert \mathbf{P}\left( AB\right) -\mathbf{P}(A)\mathbf{P}(B)\right\vert
;A\in F_{n+k}^{\infty },B\in F_{1}^{k}\right\} ,
\end{equation*}
converge to zero as $n\longrightarrow \infty $, where $F_{l}^{m}=\sigma\left\{ \zeta _{l},\zeta _{l+1},...,\zeta _{m}\right\} $ denotes the $\sigma
-$algebra generated by $\zeta _{l},\zeta _{l+1},...,\zeta _{m}$ with $l\leq m.$\\
The $\alpha $-$mixing$ condition is the most common dependence used in literature, because it is known to be fulfilled for many processes. Rosemblatt (1971) showed that a purely non-deterministic Markov chain is $\alpha$-mixing. Gorodestkii (1977) and Withers (1981) derived the conditions under which a linear process is $\alpha $-mixing. We refer the readers interested in mixing data to the book by Zhengyan and Chuanrong (2010) or to Bradley's (2007) monographs, which gives many in depth results for different types of mixing data.
\section{Assumptions and main results}
All along the paper, we denote by $C$ and $c$ any positive generic constant. For some $\tau >0$, define the set neighborhood of $m\left( x\right) ,$
$\Xi :=\left[ m\left( x\right) -\tau ,m\left( x\right) +\tau \right] ,$ and suppose that $m(x)\in \Xi \cap \left[ a_{H},\text{ }b\right] $ for any $x\in
supp(v)$, where $supp(v)$ is the support of the density $v(.)$ and $b<b_{H}.$\\
\subsection{Assumptions}
We will make use of the following assumptions needed to state our main results.\\
\textbf{Assumptions} (K)
\begin{enumerate}
\item [K1.] The kernel $K\left( .\right)$ is a bounded density function.
\item [K2.] $\int\limits_{\mathbb{R}^{d}}w_{i}K\left( w\right) dw=0,$ $\forall i=1,...,d$ and $\int\limits_{\mathbb{R}^{d}}\left\vert w_{i}\right\vert \left\vert w_{j}\right\vert K\left(w\right) dw<\infty $ $\forall i,j=1,...,d.$
\end{enumerate}
\textbf{Assumptions} (R)
\begin{enumerate}
\item[R1.] The functions $\Psi _{x}\left( u,\theta \right) $ and $\Gamma_{x}\left( u,\theta \right) :=\mathbb{E}\left[ \dfrac{\psi
_{x}^{2}\left( Y_{1}-\theta \right) }{L\left( Y_{1}\right) \overline{G}\left( Y_{1}\right) }|X_{1}=u\right] v\left( u\right) $ exists and are
continuous with respect to $u$ and $\theta $ in $\mathcal{U}\left( x\right) $ (a neighborhood of $x$) and $\Xi $ respectively, for all $x\in supp(v)$.
\item[R2.] The functions $\mathbb{E}\left[ \left\vert \psi _{x}\left(Y_{1}-\theta \right) \right\vert ^{s}|X_{1}=u\right] v\left( u\right) $ for
some $s>2$, and $\mathbf{v}\left( .\right) $ the conditional density of $X,$ w.r.t $\mathbf{P}$ exist and are uniformly bounded for $u$ in $\mathcal{U}\left( x\right) $ and $\theta \in \Xi ,$ for $x\in supp(v)$.
\item[R3.] $\mathbf{E}\left[ \left\vert \psi _{x}\left( Y_{1}-\theta\right) \psi _{x}\left( Y_{j+1}-\theta \right) \right\vert |X_{1}=u_{1},%
\text{ }X_{j+1}=u_{2}\right] \mathbf{v}_{j}\left( u_{1},u_{2}\right) $ exists and is uniformly bounded for all $u_{1}$, $u_{2}\in \mathcal{U}\left(
x\right) $ and $\theta \in \Xi ,$ for $x\in supp(v),$ where $\mathbf{v}_{j}\left(.,.\right) $ is the the joint conditional density of $X_{1}$ and $X_{j+1},$ $j\geq 1,$ w.r.t $\mathbf{P}$.
\end{enumerate}
\textbf{Assumptions} (R$^{\prime }$)\\
(R$^{\prime }$) contains assumptions R$^{\prime }$1, R$^{\prime }$2 and R$^{\prime }$3 which are analogous to R1, R2 and R3 when replacing $\psi_{x}(.)$ by its derivative function $\psi ^{\prime}_{x}(.)$.\\
\textbf{Assumptions} (E)
\begin{enumerate}
\item[E1.] The function $\Psi _{x}\left( .,\theta \right) $ is twice continuously differentiable w.r.t the first component and satisfies for all $x\in
supp(v) $ and $\theta \in \Xi $
\begin{equation*}
\underset{u\in \mathcal{U}(x)}{\sup }\text{ }\underset{\theta \in \Xi }{\sup }\left\vert\frac{\partial ^{2}\Psi_{x}}{\partial{x_{i}}\partial{x_{j}}}\left( u,\theta \right)\right\vert <\infty .
\end{equation*}
\item[E2.] The function $\psi_{x}$ is continuously differentiable, and its derivative is such that $|\psi'_{x}(t)|\geq C>0$, $\forall t\in \mathbb{R}$.
\end{enumerate}
\textbf{Assumption} M\\
The mixing coefficient satisfies $\alpha \left( k\right) =O\left(k^{-\lambda }\right)$, for some positive constant $\lambda $, where $\lambda >\max \left( 3,\dfrac{2\left( s-1\right) }{s-2}\right)$ for $s>2$,\\
\textbf{Assumptions} (B)\\
The bandwidth $h_{n}$ satisfies
\begin{enumerate}
\item[B1.] $nh_{n}^{d+4}=o(1).$
\item[B2.] There exists a sequence of positive integers $q_{n}=o\left(\sqrt{nh_{n}^{d}}\right)$
  going to infinity along with $n,$ such that $%
\underset{n\rightarrow \infty }{\lim }\left(n/h_{n}^{d}\right)^{1/2}\alpha(q_{n}) =0$.
\end{enumerate}
\textbf{Comments on the Assumptions}
\begin{enumerate}
\item Assumption K1 is quite usual in kernel estimation, and Assumption K2 is used to get the rate of the bias term.

\item Assumptions (R) deal with some regularity of the function $\Psi _{x}\left( .,.\right) $, which are used together
with M to get the expression of the Variance of $\widetilde{\Psi}_{x}\left( x,\theta \right).$ Note that condition R1 is similar to condition M9 given in Lemdani and Ould Said (2017), and conditions R2 and R3 are used in Wang and Liang (2012).
\item Assumptions (E) are used commonly in the literature of robust regression. E1 is used to get the expression of the bias, whereas condition E2 permits to use Taylor expansion.
\item Assumption B1 is needed to make the bias term negligible. B2 is needed to use a Doob's technique to get asymptotic normality and is used in Masry (2005), Wang and Liang (2012) and others. Note that condition B2 can be satisfied if we choose $q_{n}=\left[\sqrt{\frac{nh_{n}^{d}}{logn}}\right]$ and $h_{n}=O(n^{-a})$, for some $0<a<1/d,$ if the mixing coefficient satisfies $\lambda>1+(2ad/(1-ad))$.
\item Assumptions (R$^{\prime }$) are needed to get the weak convergence of $(\partial/\partial\theta)\widetilde{\Psi}_{x}\left( x,\theta \right),$ which is given by Lemma \ref{lem5}.
\end{enumerate}
\subsection{Main results}
By E2, a Taylor expansion of $\widehat{\Psi} _{x}\left( x, . \right)$ around $m(x)$, and the property
\begin{equation*}
\widehat{\Psi} _{x}\left( x,\widehat{m}\left( x\right) \right) =\Psi _{x}\left( x,m\left(
x\right) \right) =0,
\end{equation*}%
we get
\begin{equation}\label{3.1}
\Psi _{x}\left( x,m\left( x\right) \right) -\widehat{\Psi} _{x}\left( x,m\left(
x\right) \right) =\left( \widehat{m}\left( x\right) -m\left( x\right) \right) \dfrac{\partial \widehat{\Psi }}{\partial \theta }\left( x,%
\widehat{m}^{\ast }\left( x\right) \right) ,
\end{equation}
where $\widehat{m}^{\ast }\left( x\right) $ lies between $\widehat{m}\left( x\right)$ and $m(x)$. Equation (\ref{3.1}) shows that from the
behavior of $\left( \Psi _{x}-\widehat{\Psi} _{x}\right) $ as $n$ tends to infinity, it is easy to obtain asymptotic result for the sequence $%
\widehat{m}\left( .\right) ,$ (convergence, asymptotic distribution,...).\\
\subsubsection{Weak Consistency}
\begin{prop}\label{Propo1} Suppose that Conditions K1, (R) and M hold, then we have for all real $\theta\in\Xi$
  \begin{equation*}
  \widehat{\Psi }_{x}\left( x,\theta \right) \overset{\textbf{p}}{\rightarrow }\Psi_{x}\left( x,\theta \right)  \text{\ \ \textit{as}}\mathit{\ }n\rightarrow \infty.
\end{equation*}%
where $\overset{\textbf{p}}{\rightarrow }$ denotes convergence in probability w.r.t $\textbf{P}$.
\end{prop}
\vskip 4mm
\begin{theo}\label{theo1} Under Conditions of the above proposition we have
\begin{equation*}
i)\text{  }\widehat{m}\left( x\right) \overset{\textbf{p}}{\rightarrow }m\left( x\right) \text{
\ \ \textit{as}}\mathit{\ }n\rightarrow \infty. \text{\ \ \ \ \ \ \ \ \ \ \ \ \ \ \ \ \ \ \ \ \ \ \ \ \ \ \ \ \ \ \ \ \ }
\end{equation*}
Moreover, if K2 and (E) hold we get
\begin{equation*}
    ii)\text{  }\left| \widehat{m }(x)-m(x)\right|=O\left(h_{n}^{2}\right)+O_{\mathbf{p}}\left(\frac{1}{\sqrt{nh_{n}^{d}}}\right)  \text{\ \ \textit{as}}\mathit{\ }n\rightarrow \infty.
\end{equation*}
\end{theo}
\subsubsection{Asymptotic Normality}
\begin{prop}\label{Propo2}  Suppose that Conditions (K), (R), E1, M and (B)  hold, then we have
\begin{equation*}
\sqrt{nh_{n}^{d}}(\widehat{\Psi }_{x}\left( x,\theta \right) -\Psi
_{x}\left( x,\theta \right) )\overset{D}{\rightarrow }\mathcal{N}\left( 0,\sigma
_{0}^{2}\left( x,\theta \right) \right) \text{ \ \ \textit{as}}\mathit{\ }%
n\rightarrow \infty ,
\end{equation*}%
\textit{where }%
\begin{equation}\label{3.2}
\sigma _{0}^{2}\left( x,\theta \right) :=\mu\Gamma _{x}\left( x,\theta \right)
\kappa,
\end{equation}
$\Gamma _{x}\left( .,\theta \right)$ is defined in condition R1, $\kappa:=\int\limits_{\mathbb{R}^{d}}K^{2}\left( w\right) dw$, and $\overset{D}{\rightarrow }$ denotes the convergence in distribution.
\end{prop}
\vskip 4mm
\begin{theo}\label{theo2} Under Conditions (K), (R), (E), (R$^{\prime }$), M and (B) we have
\begin{equation*}
\sqrt{nh_{n}^{d}}\left( \widehat{m}\left( x\right) -m\left( x\right) \right)
\overset{D}{\rightarrow }\mathcal{N}\left( 0,\sigma ^{2}\left( x\right) \right) \text{
\ \ \textit{as}}\mathit{\ }n\rightarrow \infty ,
\end{equation*}%
\textit{where }%
\begin{equation*}
\sigma ^{2}\left( x\right) :=\dfrac{\sigma _{0}^{2}\left( x,m(x)\right) }{%
\left( \dfrac{\partial \Psi _{x}\left( x,m(x)\right) }{\partial \theta }%
\right) ^{2}},
\end{equation*}
and $\sigma _{0}^{2}\left( x,.\right) $ is defined in (\ref{3.2}).
\end{theo}
\vskip 4mm
\begin{rmk}\label{rem1} \textbf{Special cases:}
\begin{enumerate}
  \item \textbf{Censored data:} in absence of truncation $(T=0)$, our results extend the results given by Lemdani and Ould said (2017) in the case of independent data, and the asymptotic variance becomes
      \begin{equation*}
        \sigma ^{2}\left( x\right) :=\dfrac{\mathbb{E}\left[ \psi_{x}^{2}\left(Y-m(x)\right)\overline{G}^{-1}\left( Y\right) |X=x\right]v(x)\kappa }{%
\left( \dfrac{\partial \Psi _{x}\left( x,m(x)\right) }{\partial \theta }\right) ^{2}},
      \end{equation*}
      which is the same as that obtained in their Theorem 3.2.
  \item \textbf{Truncated data:} in absence of censoring $(W=\infty)$, the asymptotic variance becomes
  \begin{equation*}
        \sigma ^{2}\left( x\right) :=\dfrac{\mu\mathbb{E}\left[ \psi_{x}^{2}\left(Y-m(x)\right){L}^{-1}\left( Y\right) |X=x\right]v(x) \kappa}{%
\left( \dfrac{\partial \Psi _{x}\left( x,m(x)\right) }{\partial \theta }\right) ^{2}},
      \end{equation*}
which is compared to that given by Wang and Liang (2012), where $\mathbb{E}\left[\dfrac{d\psi _{x}\left( Y-m(x)\right) }{d\theta }|X=x\right]
$ represent $\lambda _{1}\left( x\right) $ defined in their Theorem 2.3. Note that our theorem holds with a relaxed condition on the mixing coefficient. (Theorem 2.3 of Wang and Liang (2012) holds with $\lambda >\dfrac{s(s+2)}{2(s-2)}$, which increases along with $s$, while our theorem states with $\lambda >\max \left( 3,\dfrac{2(s-1)}{s-2}\right), $ for $s>2,$ which reduces to $\lambda >3$ for $s>4$, where $s=\infty$ correspond to the case of bounded objective function $\psi(.)$).
\end{enumerate}
\end{rmk}
\subsubsection{Application to Confidence Interval}
A direct application of Theorem \ref{theo2} allows to get a confidence interval for the true value of $m(x)$, given $X = x$. Notice that Theorem \ref{theo2} is useless in practice since many quantities in the variance are unknown. Therefore replacing the two functions $\dfrac{\partial
\Psi _{x}\left( x,m(x)\right) }{\partial \theta }$ and $\Gamma _{x}\left(x,m\left( x\right) \right) $ by their unbiased estimates defined
respectively by
\begin{eqnarray}
\dfrac{\partial \widetilde{\Psi }_{x}\left( x,m(x)\right) }{\partial \theta }
&:=&\frac{\mu }{nh_{n}^{d}}\sum\limits_{i=1}^{n}K\left( \frac{x-X_{i}}{h_{n}}%
\right) \frac{\delta _{i}\frac{\partial \psi _{x}\left( Z_{i}-m\left(
x\right) \right) }{\partial \theta }}{L\left( Z_{i}\right) \overline{G}%
\left( Z_{i}\right) }, \label{3.4}\\
\widetilde{\Gamma }_{x}\left( x,m\left( x\right) \right) &:=&\frac{\mu}{%
nh_{n}^{d}}\sum\limits_{i=1}^{n}K\left( \frac{x-X_{i}}{h_{n}}\right) \frac{%
\delta _{i}\psi _{x}^{2}\left( Z_{i}-m\left( x\right) \right) }{L^{2}\left(
Z_{i}\right) \overline{G}^{2}\left( Z_{i}\right) },  \label{3.4bis}
\end{eqnarray}
and using a plug-in method we obtain a convergent estimate $\widehat{\Gamma }_{x}\left( x,\widehat{m}\left( x\right) \right) $ and $\dfrac{\partial
\widehat{\Psi }_{x}\left( x,\widehat{m}(x)\right) }{\partial \theta }$ by substituting  $\mu,$\, $\overline{G},$\, $L,$ and $m(x)$ in (\ref{3.4}) and (\ref{3.4bis}) by their estimates $\mu _{n},$\ $\overline{G}_{n},$\, $L_{n},$ and $\widehat{m}(x)$ given by (\ref{2.2}), (\ref{G}), (\ref{L}) and solution of (\ref{2.6}) respectively, which permit to obtain a convergent estimate $\widehat{\sigma }$\ of $\sigma $\ defined by
\begin{equation}\label{3.5}
\widehat{\sigma }^{2}\left( x\right) :=\frac{\mu_{n} \widehat{\Gamma }_{x}\left( x,%
\widehat{m}\left( x\right) \right)\kappa }{\left( \dfrac{\partial \widehat{\Psi }%
_{x}\left( x,\widehat{m}(x)\right) }{\partial \theta }\right) ^{2}}%
,
\end{equation}
if the denominator does not vanish. Then, from Theorem \ref{theo2} we get the following
\begin{corol}\label{corol1} Under the assumptions of Theorem \ref{theo2}, we have
for $n\rightarrow \infty $
\begin{equation*}
\frac{\sqrt{nh_{n}^{d}}}{\widehat{\sigma }\left( x\right) }\left( \widehat{m}%
\left( x\right) -m\left( x\right) \right) \overset{D}{\rightarrow }\mathcal{N}\left(
0,1\right) .
\end{equation*}
\end{corol}
From this corollary, we get for each fixed $\eta \in \left( 0,1\right) ,$ the following $\left( 1-\eta \right) $ asymptotic level confidence interval
\begin{equation*}
\left[ \widehat{m}\left( x\right) -t_{1-\eta /2}\times \frac{\widehat{\sigma
}\left( x\right) }{\sqrt{nh_{n}^{d}}}\text{ },\text{ }\widehat{m}\left(
x\right) +t_{1-\eta /2}\times \frac{\widehat{\sigma }\left( x\right) }{\sqrt{%
nh_{n}^{d}}}\right] ,
\end{equation*}
where $t_{1-\eta /2}$ denotes the $\left(1-\eta /2\right)$ quantile of the standard normal distribution.
\subsubsection{Comeback to Classical Regression}
As a direct consequence of Theorem \ref{theo1} and Theorem \ref{theo2}, we get the weak consistency and asymptotic normality of the estimation of the regression curve $\mathbb{E[}Y|X=.],$ stated as solution of equation (\ref{2.6}) with the identity objective function $\psi_{x}(u)=u$ given by
\begin{equation}\label{3.6}
\widehat{m}_{n}(x):=\frac{\sum\limits_{i=1}^{n}K\left( \frac{x-X_{i}}{h_{n}}%
\right) \delta _{i}Z_{i}L_{n}^{-1}\left( Z_{i}\right) \overline{G}%
_{n}^{-1}\left( Z_{i}\right) }{\sum\limits_{i=1}^{n}K\left( \frac{x-X_{i}}{%
h_{n}}\right) \delta _{i}L_{n}^{-1}\left( Z_{i}\right) \overline{G}%
_{n}^{-1}\left( Z_{i}\right) },
\end{equation}
which yield to the following results
\begin{corol}\label{corol2} Under the assumptions of Theorem \ref{theo1}, we have
\begin{equation*}
\widehat{m}_{n}\left( x\right) \overset{\textbf{p}}{\rightarrow }m\left( x\right) \text{
\ \ \textit{as}}\mathit{\ }n\rightarrow \infty .
\end{equation*}
\end{corol}
\begin{corol}\label{corol3} Under the assumptions of Theorem \ref{theo2}, we have
\begin{equation*}
\frac{\sqrt{nh_{n}^{d}}}{\sigma _{class}\left( x\right) }\left( \widehat{m}%
_{n}\left( x\right) -m\left( x\right) \right) \overset{D}{\rightarrow }%
\mathcal{N}\left( 0,1\right)\text{\ \ \textit{as}}\mathit{\ }n\rightarrow \infty,
\end{equation*}
where%
\begin{equation}\label{3.7}
\sigma _{class}^{2}\left( x\right) :=\mu\mathbb{E}\left[ \dfrac{ \left(
Y-m(x)\right) ^{2}}{L\left( Y\right) \overline{G}\left( Y\right) }|X=x%
\right] v^{-1}\left( x\right) \kappa.
\end{equation}
\end{corol}
\begin{rmk}\label{rem3}
Note that the variance $\sigma _{class}^{2}\left( x\right) $ given in (\ref{3.7}) reduces to the definitions below in the case of:
\begin{enumerate}
  \item \textbf{Truncated data:} In absence of censoring $\left(\overline{G}\left( .\right) =1\right),$ we get
  \begin{equation*}
  \sigma _{class, trun}^{2}\left( x\right) :=\frac{\mu \mathbb{E}\left[\left(
Y-m(x)\right) ^{2}L^{-1}\left( Y\right) |X=x\right]\kappa }{v(x)}.
  \end{equation*}
  This is what it was obtained by Liang (2011).
  \item \textbf{Censored data:} In absence of truncation $\left( L\left( .\right) =1\right)$,  we get
\begin{equation*}
\sigma _{class, cens}^{2}\left( x\right) :=\frac{\mathbb{E}\left[ \left(
Y-m(x)\right) ^{2}\overline{G}^{-1}\left( Y\right) |X=x\right]\kappa }{v(x)}.
\end{equation*}
Which is different from that obtained by Guessoum and Ould-Said (2012), for $a_{F}<a_{G}$, and given by
\begin{equation*}
\sigma _{1}^{2}\left( x\right) =\frac{\left(\mathbb{E}\left[ Y^{2}\overline{G}%
^{-1}\left( Y\right) |X=x\right] -m^{2}(x)\right)\kappa}{v\left( x\right) }.
\end{equation*}
\end{enumerate}
\end{rmk}
Note that for every $x\in supp(v)$, we have $\sigma _{class, cens}^{2}\left( x\right) \leq\sigma _{1}^{2}\left( x\right)$, as shown in the following Proposition.
\begin{prop}\label{Propo3} For every $x\in supp(v)$ such that $m(x)$ exists, for $a_{F}<a_{G}$ we have
\begin{equation*}
\sigma _{1}^{2}\left( x\right) =\sigma _{class, cens}^{2}\left( x\right)+%
\frac{A\left( x\right) }{v\left( x\right) }\kappa,
\end{equation*}
where
\begin{eqnarray*}
A\left( x\right) &:=&\mathbb{E}\left[ (2Y-m(x))m(x)\overline{G}^{-1}\left(Y\right) |X=x\right] -m^{2}(x)\geq0.
\end{eqnarray*}
The equality $\sigma _{class, cens}^{2}\left( x\right)=\sigma _{1}^{2}\left( x\right)$ is achieved only if $m(x)=0$.
\end{prop}
{\bf Proof:}  It suffices to observe that
\begin{eqnarray*}
\mathbb{E}\left[ Y^{2}\overline{G}^{-1}\left( Y\right) |X=x\right] -m^{2}(x)
&=&\mathbb{E}\left[ \left( Y-m(x)\right) ^{2}\overline{G}^{-1}\left(
Y\right) |X=x\right] +\mathbb{E}\left[ 2Ym(x)\overline{G}^{-1}\left(
Y\right) |X=x\right] \\
&&-\mathbb{E}\left[ m^{2}(x)\overline{G}^{-1}\left( Y\right) |X=x\right]
-m^{2}(x) \\
&=:&\mathbb{E}\left[ \left( Y-m(x)\right) ^{2}\overline{G}^{-1}\left(
Y\right) |X=x\right] +A\left( x\right),
\end{eqnarray*}
where
\begin{eqnarray*}
A\left( x\right) &\geq&\overline{G}^{-1}\left( a_{F}\right) \left[ m(x)\mathbb{E%
}\left[ 2Y|X=x\right] -m^{2}(x)\right] -m^{2}(x) \\
&\geq&m^{2}(x)\left( \overline{G}^{-1}\left( a_{F}\right) -1\right) \geq0.
\end{eqnarray*}
$\hfill\Box$\\
Then, the estimator obtained from (\ref{3.6}) in absence of truncation and given by
\begin{equation*}
m_{adjNW}(x):=\frac{\sum\limits_{i=1}^{n}K\left( \frac{x-X_{i}}{h_{n}}%
\right) \delta _{i}Z_{i}\overline{G}_{n}^{-1}\left( Z_{i}\right) }{%
\sum\limits_{i=1}^{n}K\left( \frac{x-X_{i}}{h_{n}}\right) \delta _{i}%
\overline{G}_{n}^{-1}\left( Z_{i}\right) },
\end{equation*}
has a smaller asymptotical variance than that of the Carbonez estimator given by
\begin{equation*}
m_{n}(x):=\frac{\sum\limits_{i=1}^{n}K\left( \frac{x-X_{i}}{h_{n}}\right)
\delta _{i}Z_{i}\overline{G}_{n}^{-1}\left( Z_{i}\right) }{\sum\limits_{i=1}^{n}K\left(
\frac{x-X_{i}}{h_{n}}\right) }.
\end{equation*}
 \begin{rmk}
The result given in Proposition \ref{Propo3} is very important in the context of censored data, since it allows us to deduce other conditional functional estimators (such as, conditional quantile, conditional density, conditional mode, ...) which are more consistent than those existing in the literature.
\end{rmk}

\section{Simulation study}
In this section, we investigate the goodness-of-fit to the normal distribution expected from our theoretical results in section 3 on the
M-Estimator $\widehat{m}\left( x\right),$ as well as the curve of 95\% confidence intervals of $m\left( x\right).$ We do this by considering some
fixed sample-size which is the case in practice, for different percentage of truncation and censoring. For that we consider the strong mixing
bidimensional processes $\left( X_{t},Y_{t}\right) $ generated by

\[
\left\{
\begin{array}{c}
X_{1}=0.5e_{1} \\
X_{t+1}=\rho X_{t}+0.5e_{t+1} \\
Y_{t}=m(X_{t})+\epsilon _{t}\text{ \ \ \ \ }t\geq 1,%
\end{array}%
\right.
\]
where $e_{t},$ $t\geq 1,$ is independently distributed from standard Normal random variable, and for each $t,$ $e_{t}$ is independent from $X_{t}$ , $\epsilon _{t}$ is a white noise sequences with gaussian distribution $\mathcal{N}(0,\sigma ^{2})$, we take $\sigma =0.1,$ and $0<\rho <1$ is some constant, which is chosen to control the dependency of the observations, we take $\rho=0.9$ (which is the case of high dependency).\\
The censoring and truncation distributions are simulated from Exponential $ \mathcal{E}(a_{0})$ and Normal $\mathcal{N}(u_{0},2)$ distributions, respectively, where the parameters $a_{0}$ and $u_{0}$ are chosen to get different rate of censoring and truncation. The quadruples $(X,Y,T,W)$ were drawn until $n$ of them satisfiy the condition $T\leq Z=min(Y,W)$. In this way, censored and truncated samples
$\left( X_{i},Z_{i},T_{i},\delta _{i}\right),$ $i=1,...,n$ were obtained. We use the standard Gaussian kernel, and the objective function $\psi \left(.\right) $ defined by $\psi \left( u\right) =\dfrac{u}{\sqrt{1+u^{2}}}.$\\

We first examine the shape of the estimated density and the QQ-plots of $\widehat{m}\left( x\right) $ (with normalized deviation), compared to the theoretical normal-distribution at $x=0$. We draw $B$ independent $n$ samples, and we calculate $B=200$ estimations, $%
\widehat{m}_{i}\left( 0\right) $, $1\leq i\leq 200$ of $m\left( 0\right),$ we take the linear model $m(x)=2x.$ We derive an estimator $\widehat{%
\sigma }\left( 0\right) $ of $\sigma \left( 0\right) $ as defined in (\ref{3.5}) and the normalized deviation between this estimate and the theoretical
regression function as detailed in corollary \ref{corol1}, that is,
\[
m_{n}:=\frac{\sqrt{nh_{n}}}{\widehat{\sigma }\left( 0\right) }\left(
\widehat{m}\left( 0\right) -m\left( 0\right) \right) =\frac{\sqrt{nh_{n}}}{%
\widehat{\sigma }\left( 0\right) }\widehat{m}\left( 0\right).
\]

The bandwidth $h_{n}$ is chosen as the minimizer of the Least Squares Cross Validation (LSCV) criterion, which takes values around $1.13$. We draw, using this scheme, $B$ independent $n$-samples. In order to estimate the density function of $m_{n}$ by the kernel method, we make the classical bandwidth
choice (see, e.g., Silverman (1986), p.40) $h_{B}\approx1.6n^{-1/5}$.

We have drawn several graphs of the density of the $m_{n}$ against a standard normal distribution and the corresponding QQ-plot from different values of $n$, censoring rate $CR$ and truncation rate $TR$.\\
\begin{figure}[H]
\begin{center}
 \includegraphics[scale=0.36]{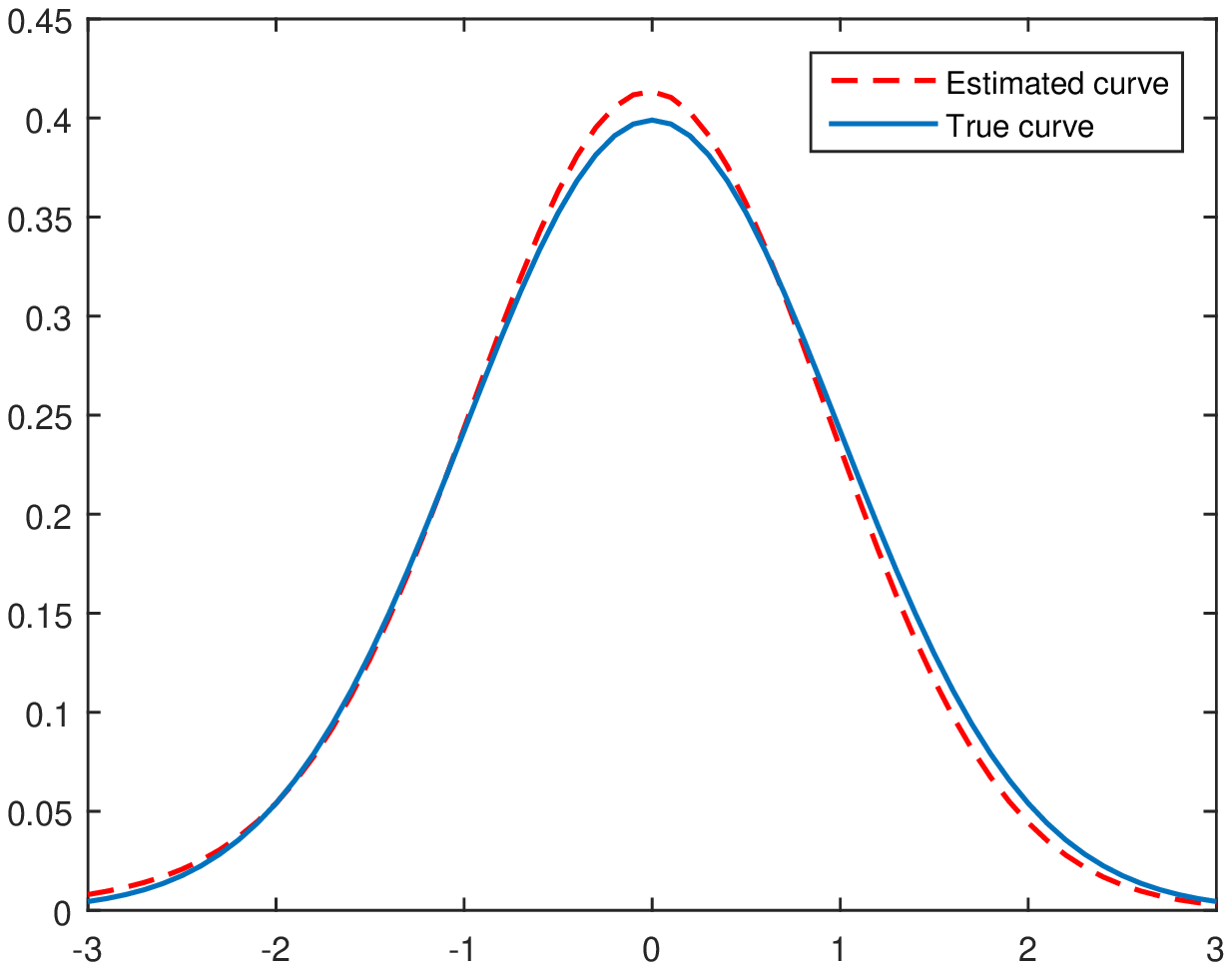}
 \includegraphics[scale=0.36]{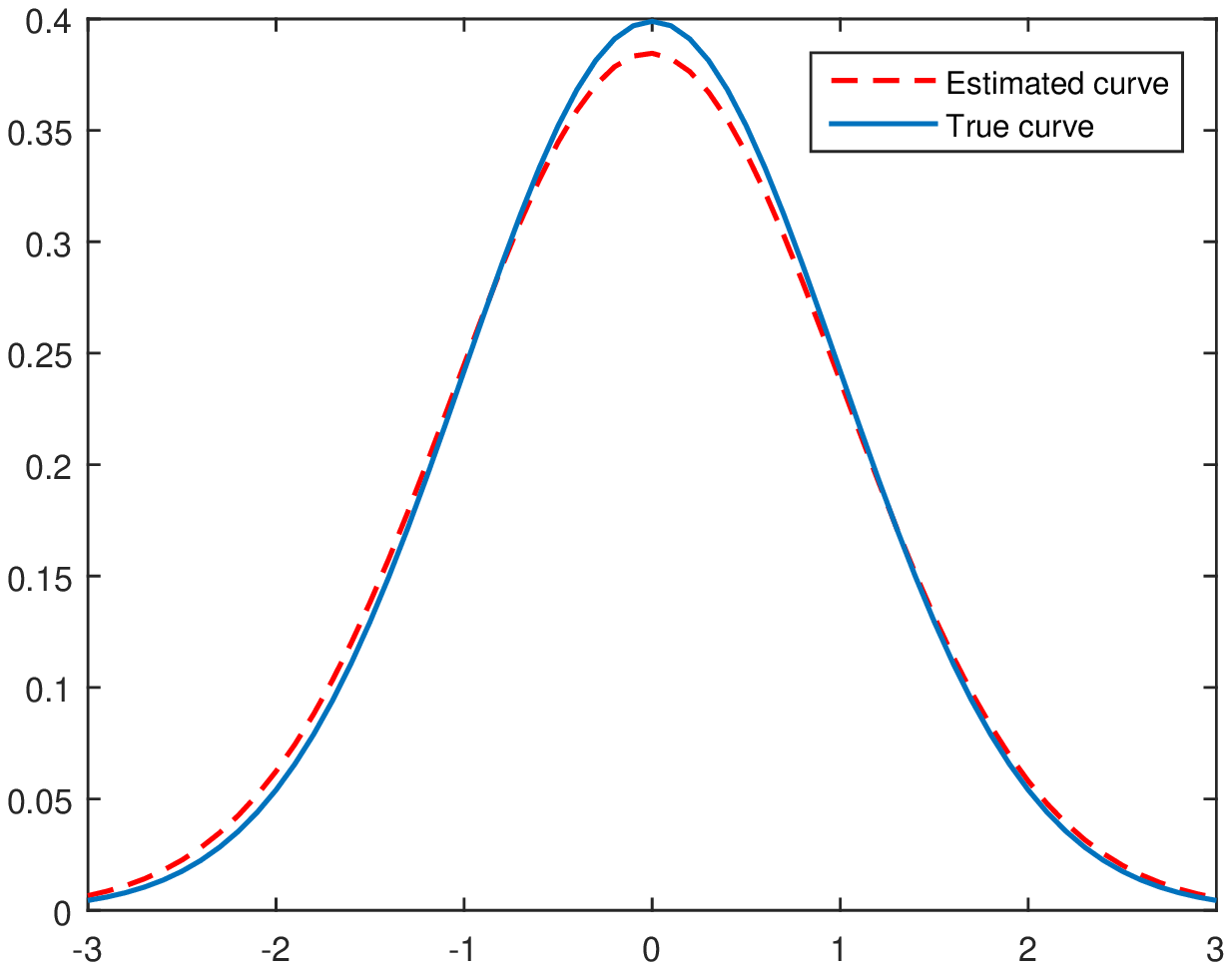}
 \includegraphics[scale=0.36]{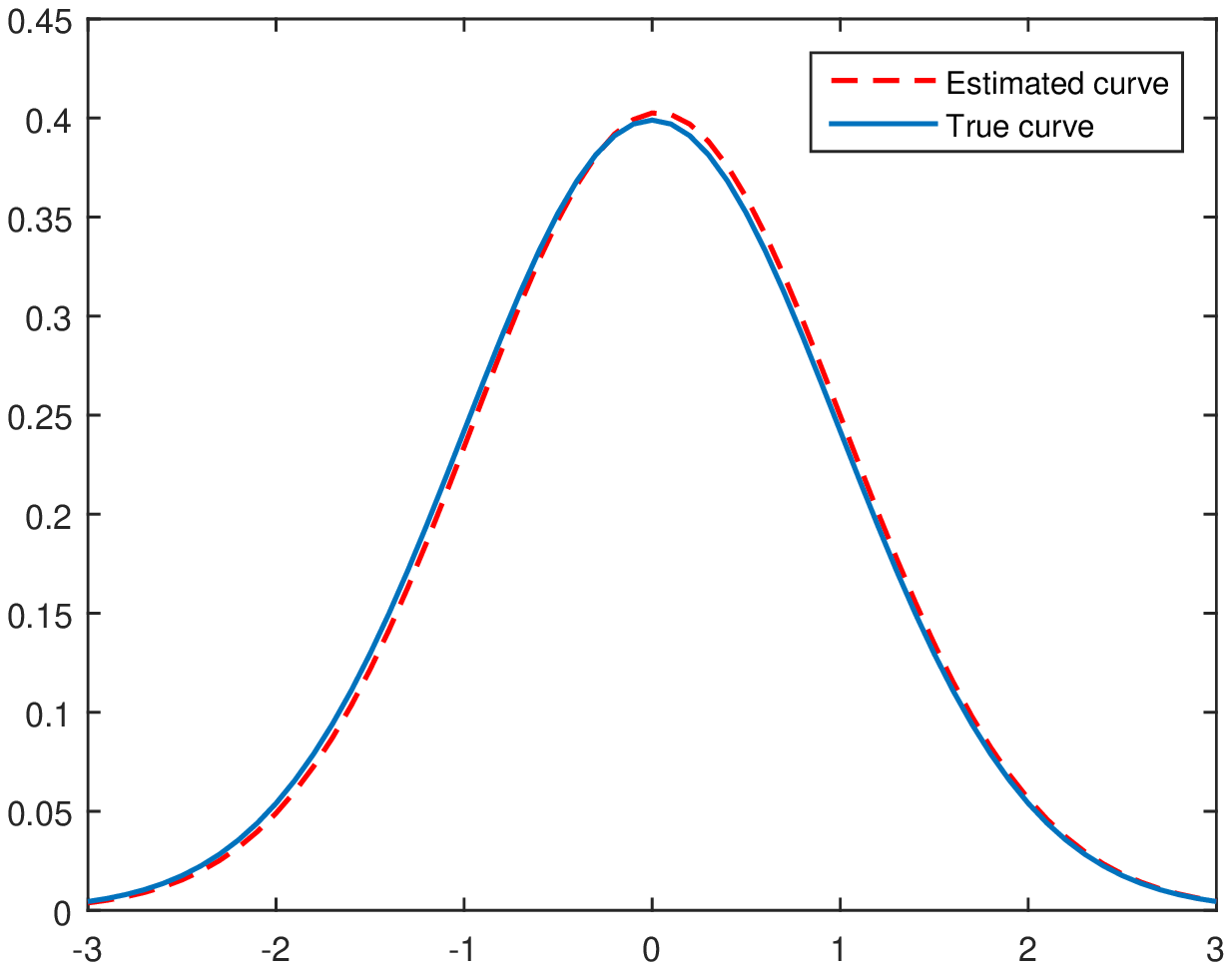}
\caption{\footnotesize{$\rho =0.9,$ $TR\approx20\%$, $ CR\approx 20\%,$ $B=200$ for $n= 50$, $100$ and $300$, respectively.}}\label{Figure5.1}
 \end{center}
\end{figure}
\begin{figure}[H]
\begin{center}
 \includegraphics[scale=0.36]{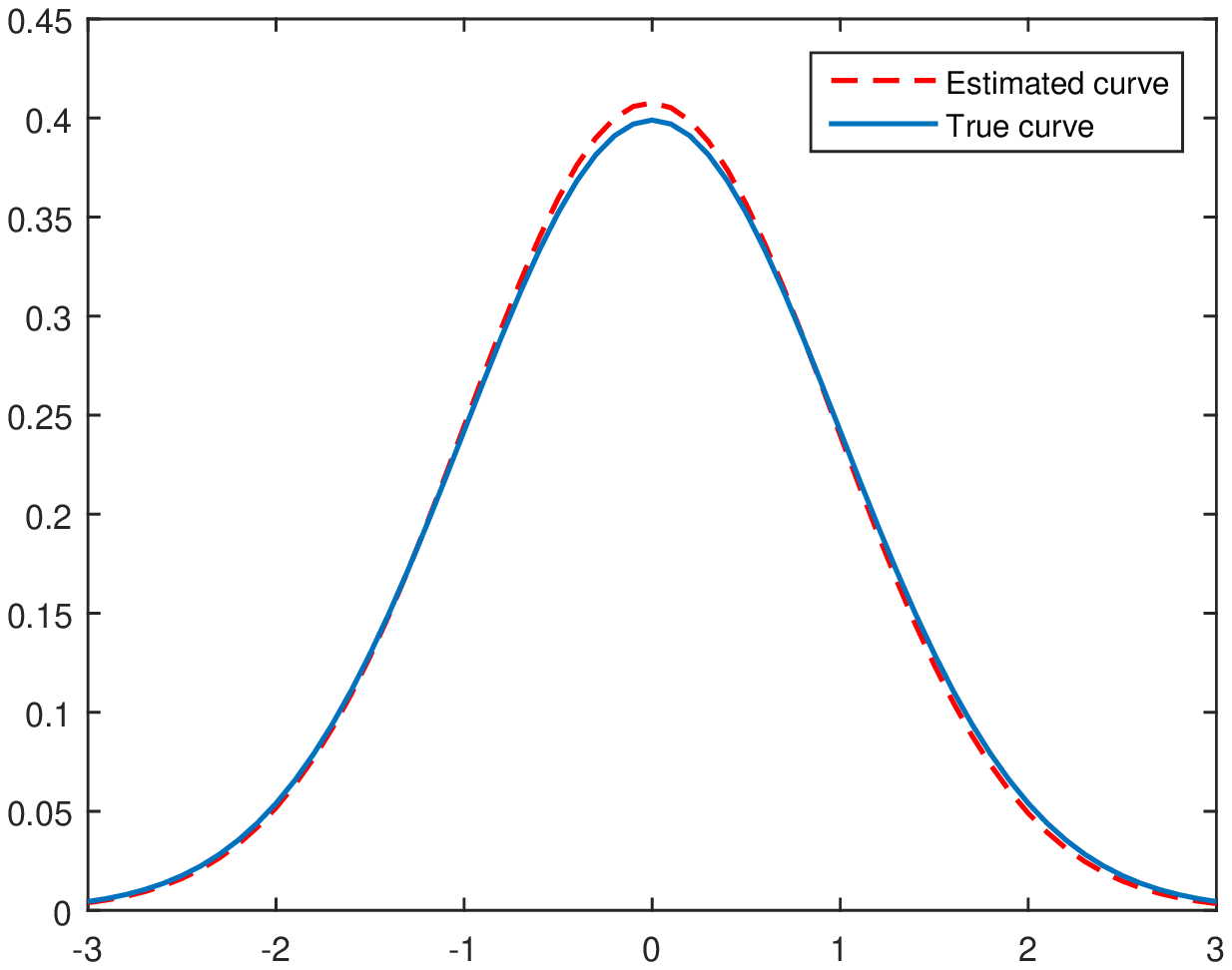}
 \includegraphics[scale=0.36]{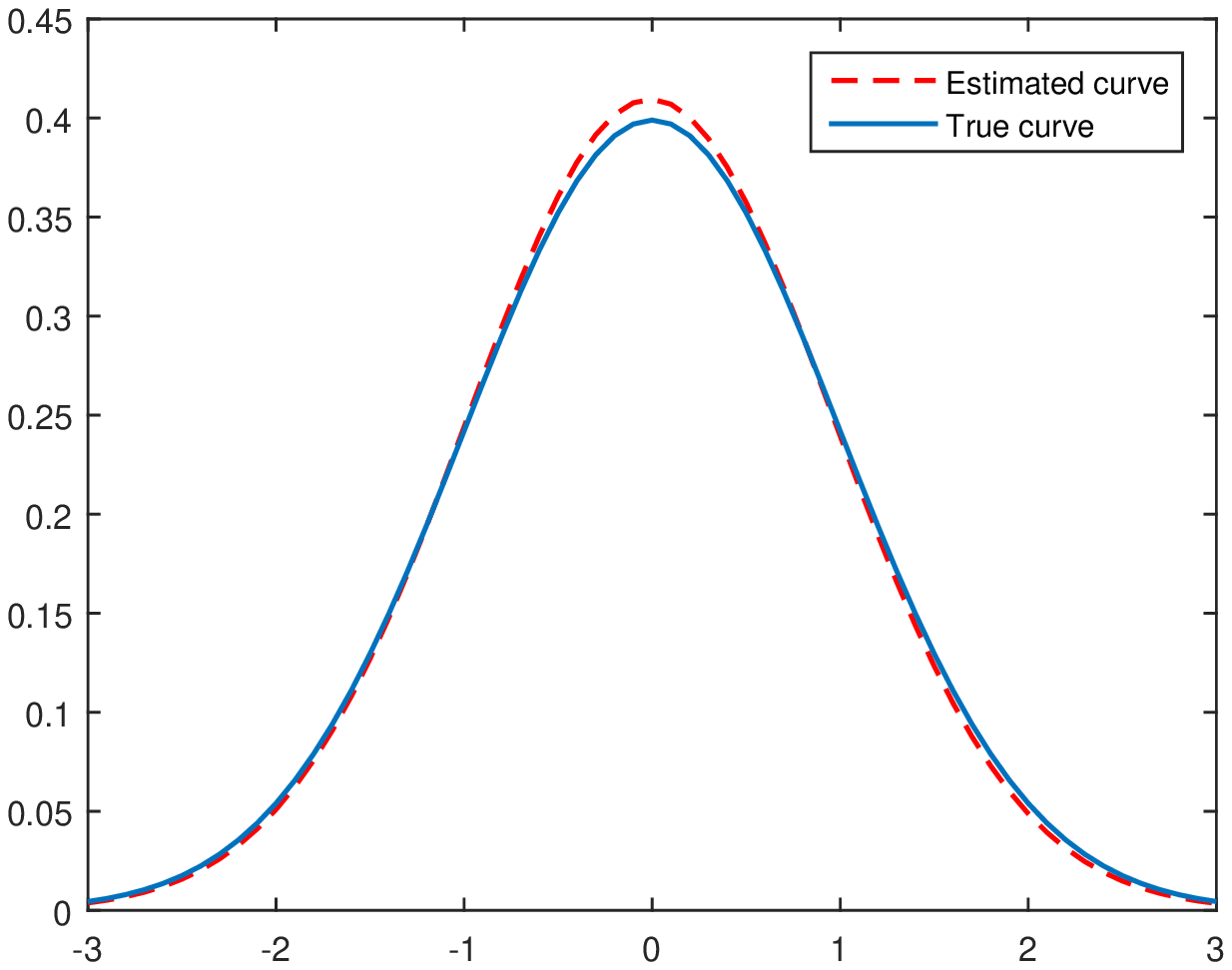}
 \includegraphics[scale=0.36]{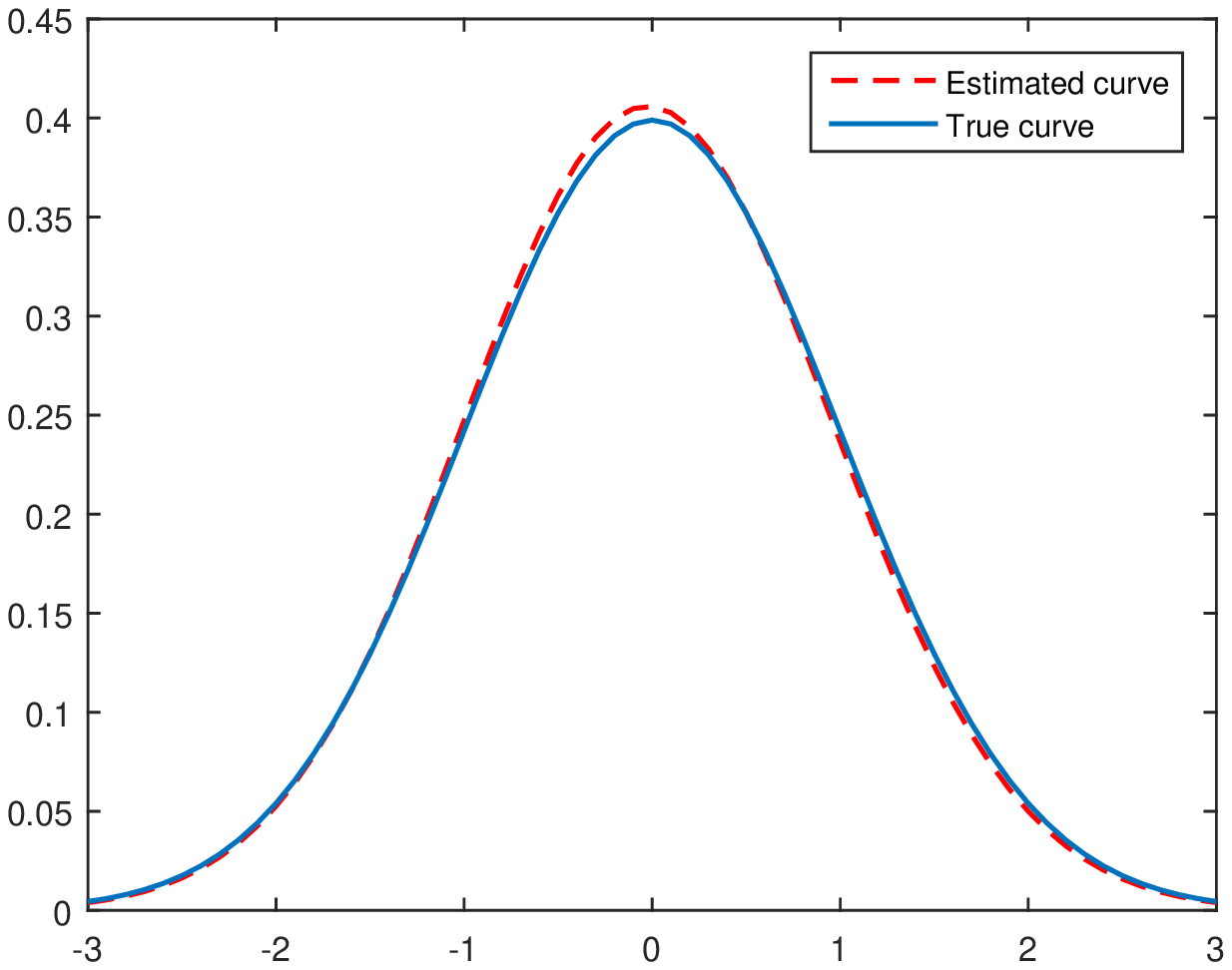}
\caption{\footnotesize{$\rho =0.9$, $TR\approx20\%$, $n= 300$, $B=200$ for $ CR\approx 0\%,$ $10\%$ and $40\%$, respectively.}}\label{Figure5.2}
 \end{center}
\end{figure}
\begin{figure}[H]
\begin{center}
 \includegraphics[scale=0.36]{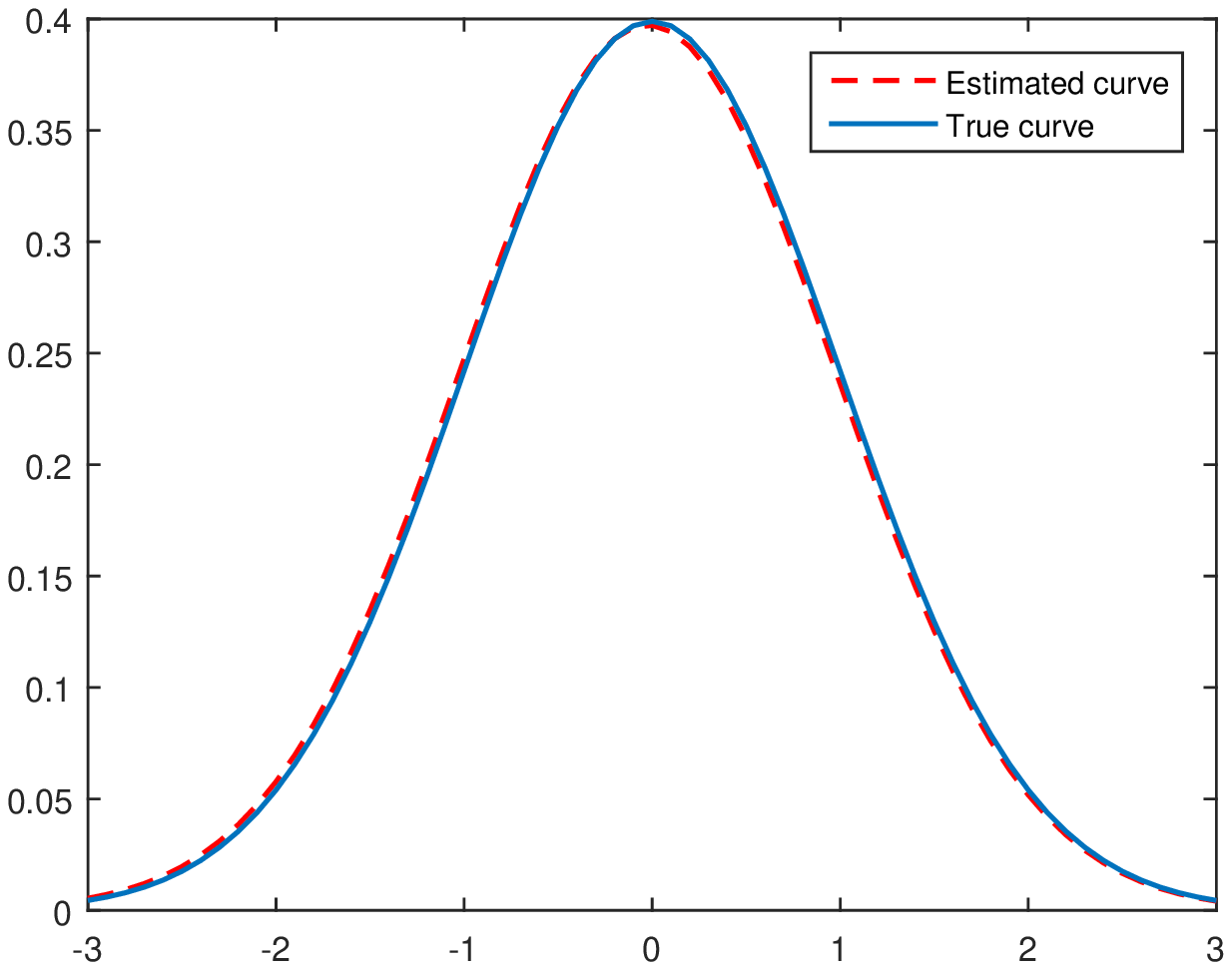}
 \includegraphics[scale=0.36]{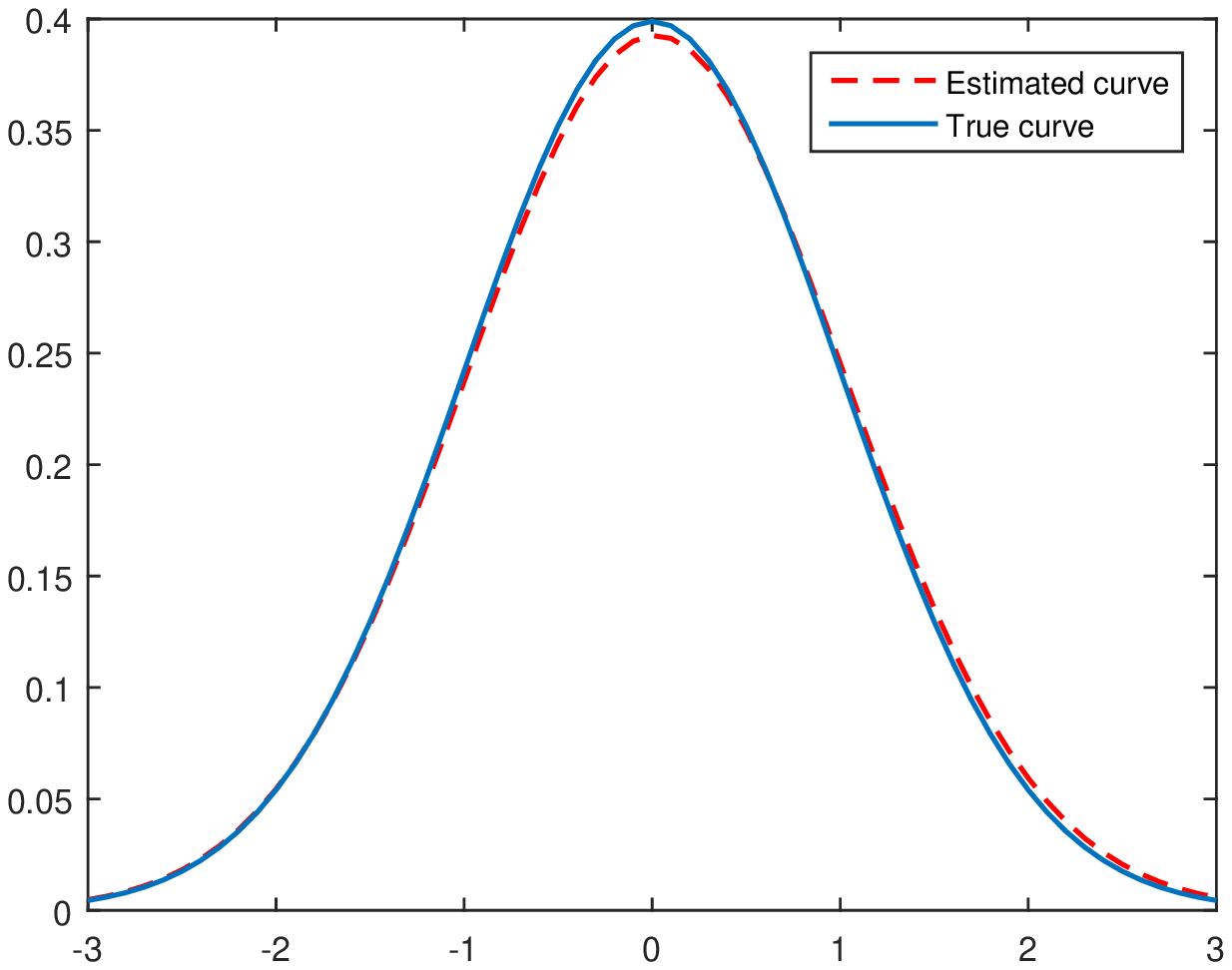}
 \includegraphics[scale=0.36]{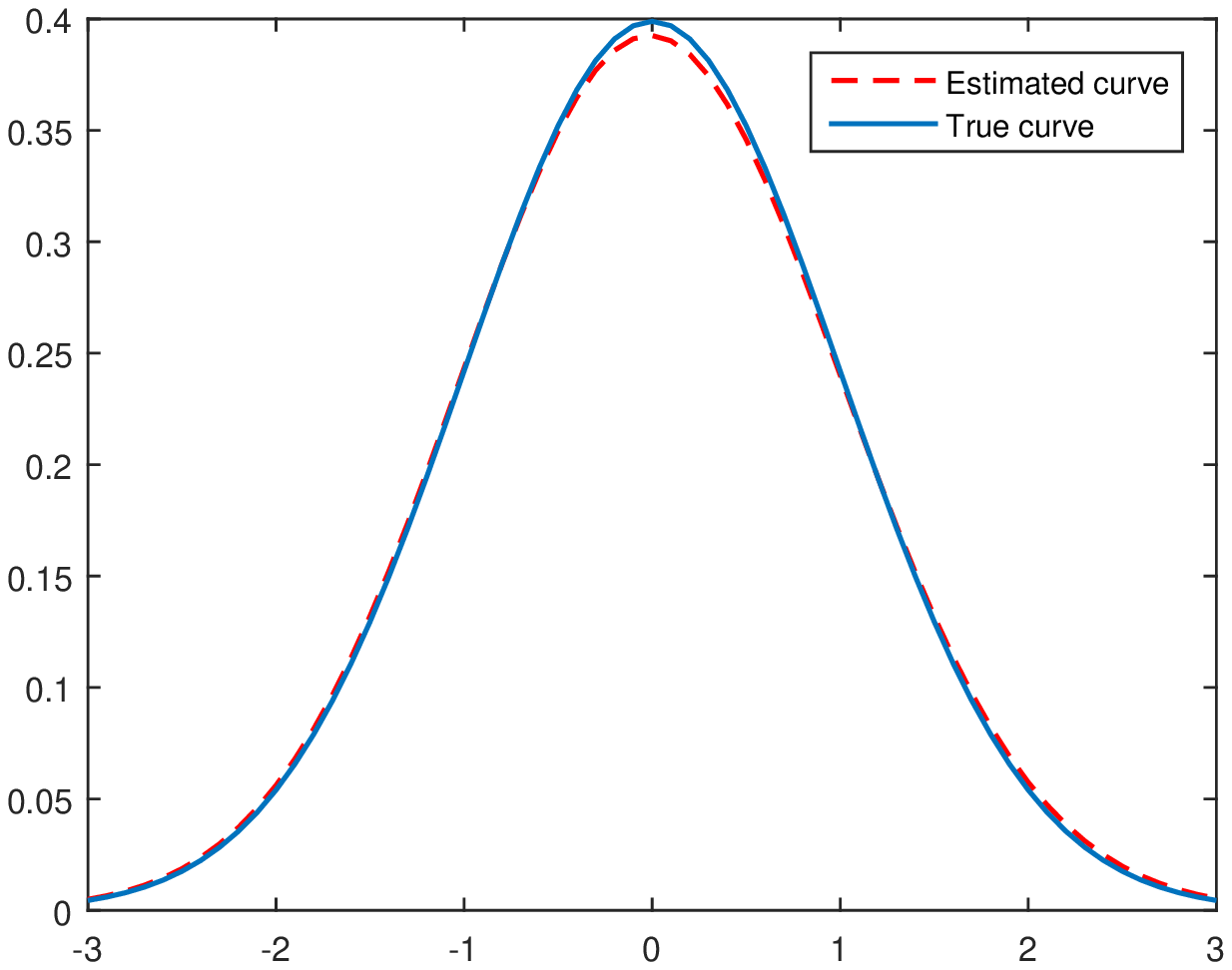}
\caption{\footnotesize{$\rho =0.9$, $CR\approx20\%$, $n= 300$, $B=200$ for $ TR\approx 0\%,$ $10\%$ and $40\%$, respectively.}}\label{Figure5.3}
 \end{center}
\end{figure}
\begin{figure}[H]
\begin{center}
 \includegraphics[scale=0.36]{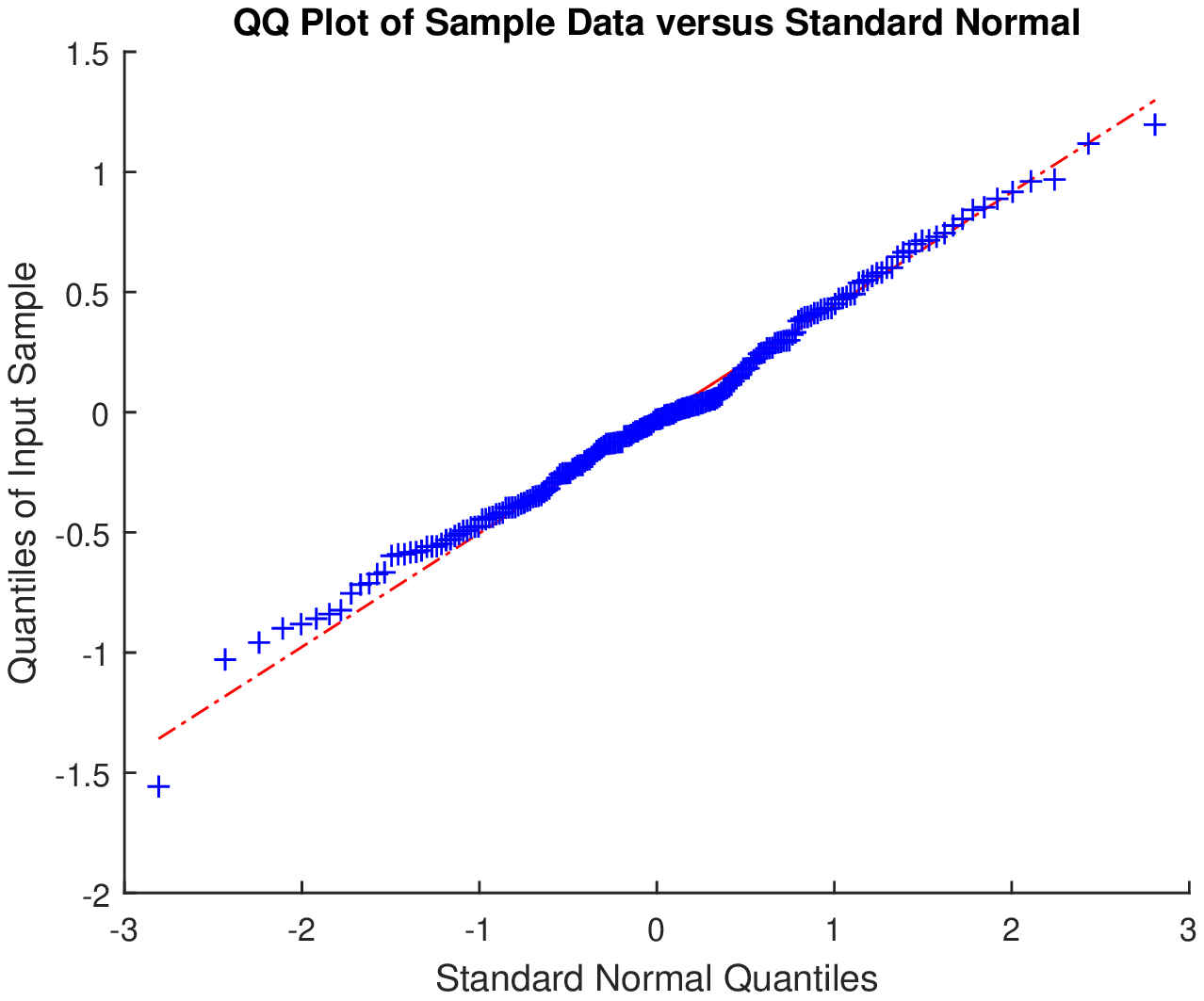}
 \includegraphics[scale=0.36]{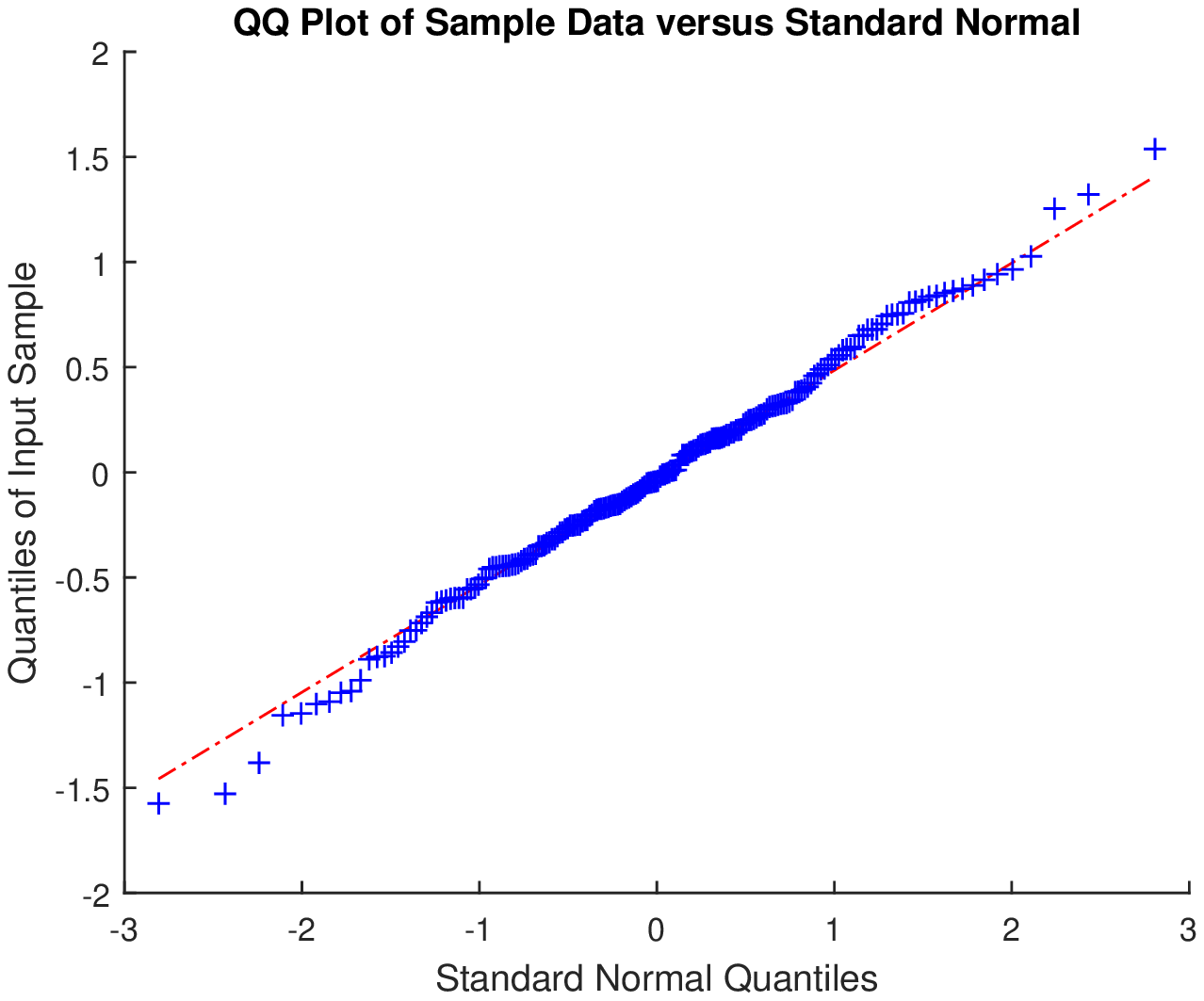}
 \includegraphics[scale=0.36]{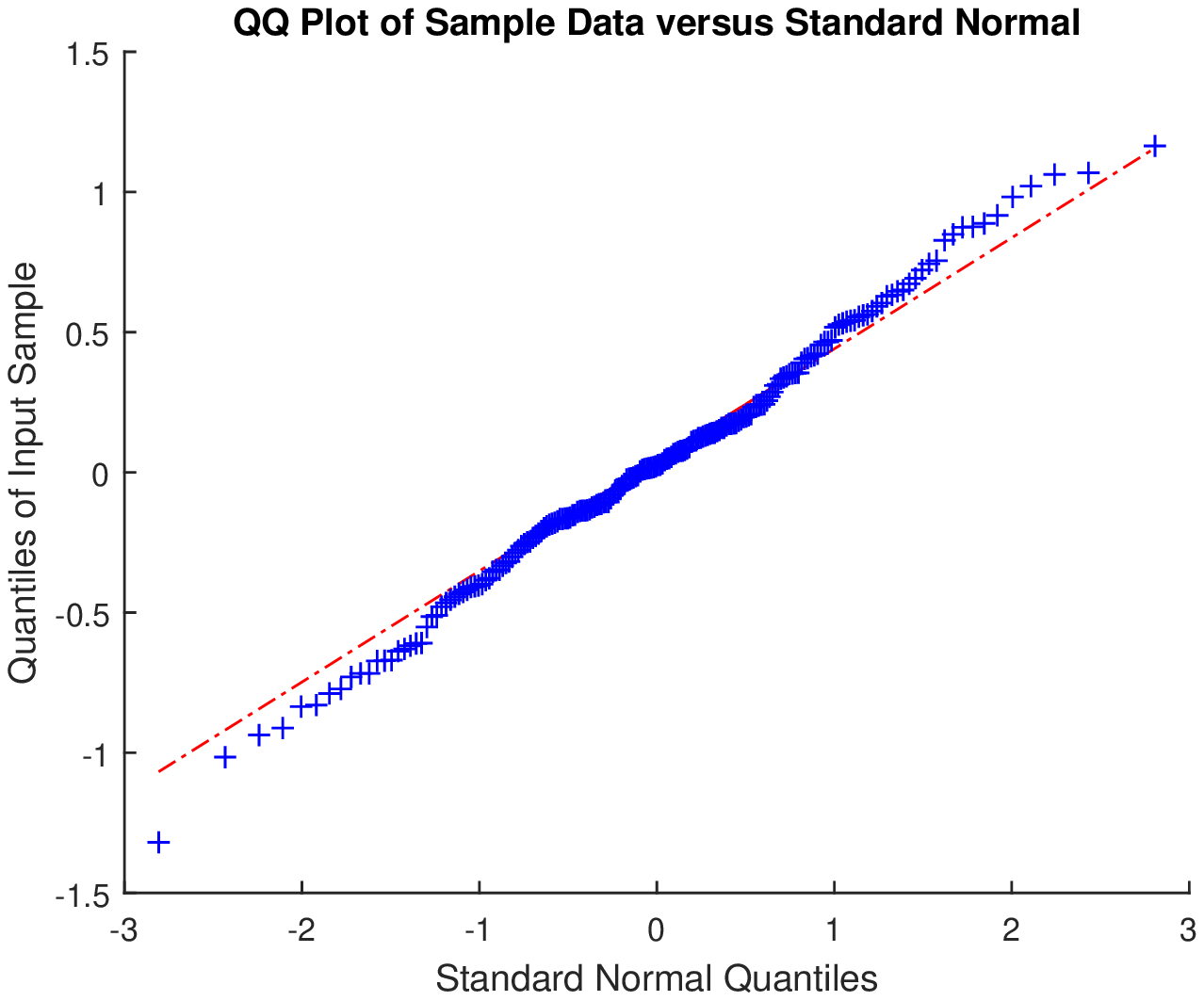}
\caption{\footnotesize{$\rho =0.9,$ $TR\approx20\%$, $ CR\approx 20\%,$ $B=200$ for $n= 50$, $100$ and $300$, respectively.}}\label{Figure5.4}
 \end{center}
\end{figure}
\begin{figure}[H]
\begin{center}
 \includegraphics[scale=0.36]{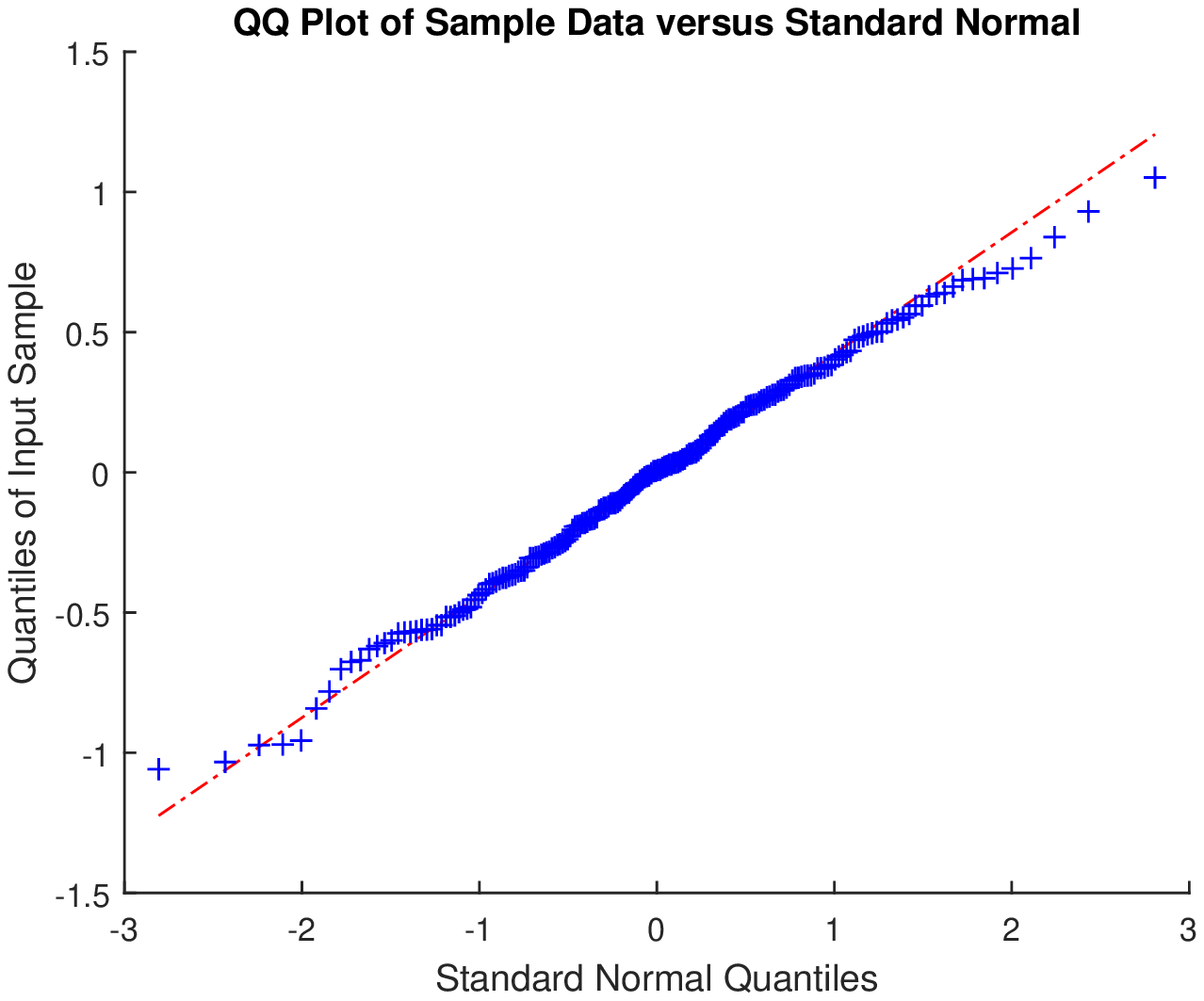}
 \includegraphics[scale=0.36]{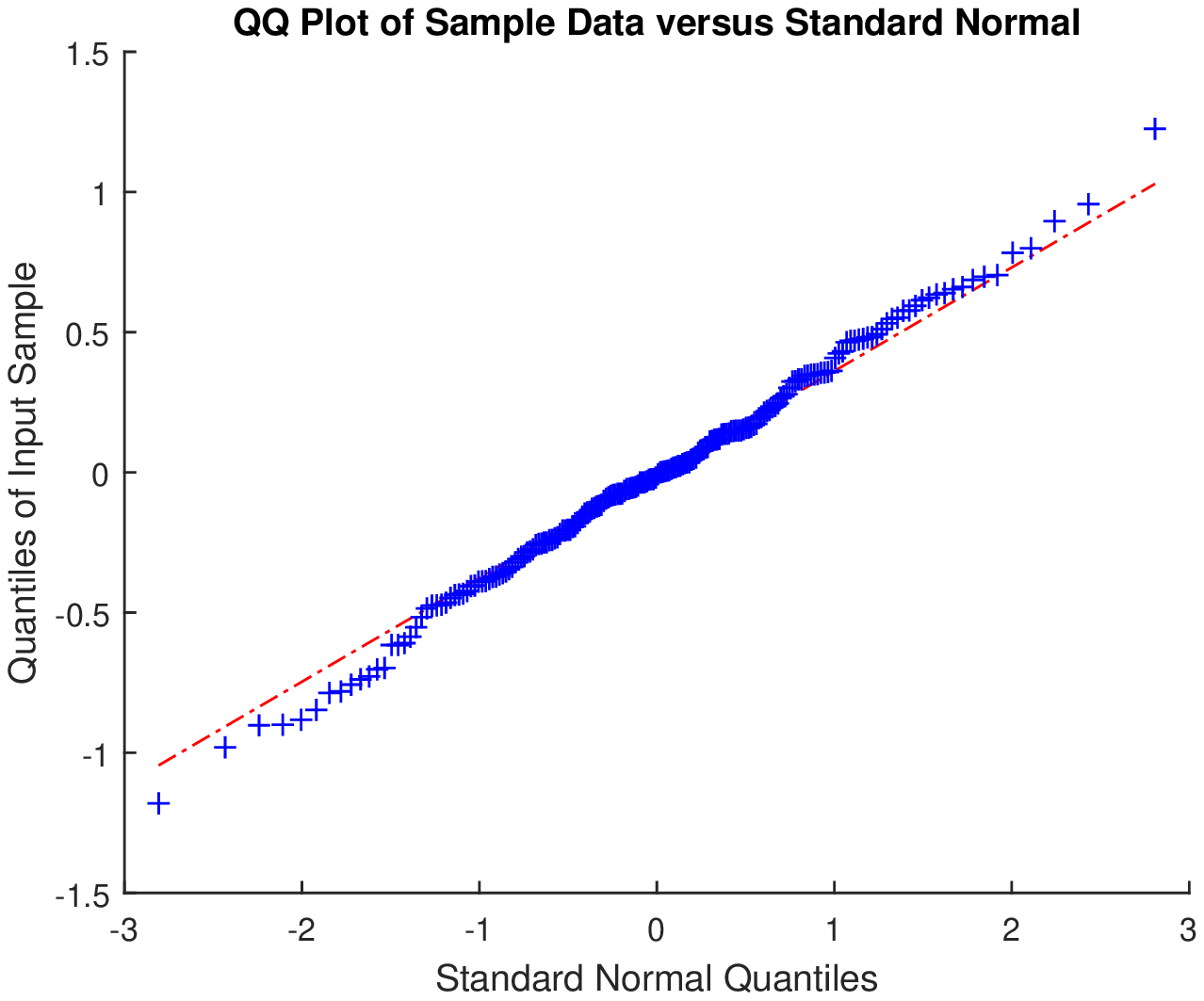}
 \includegraphics[scale=0.36]{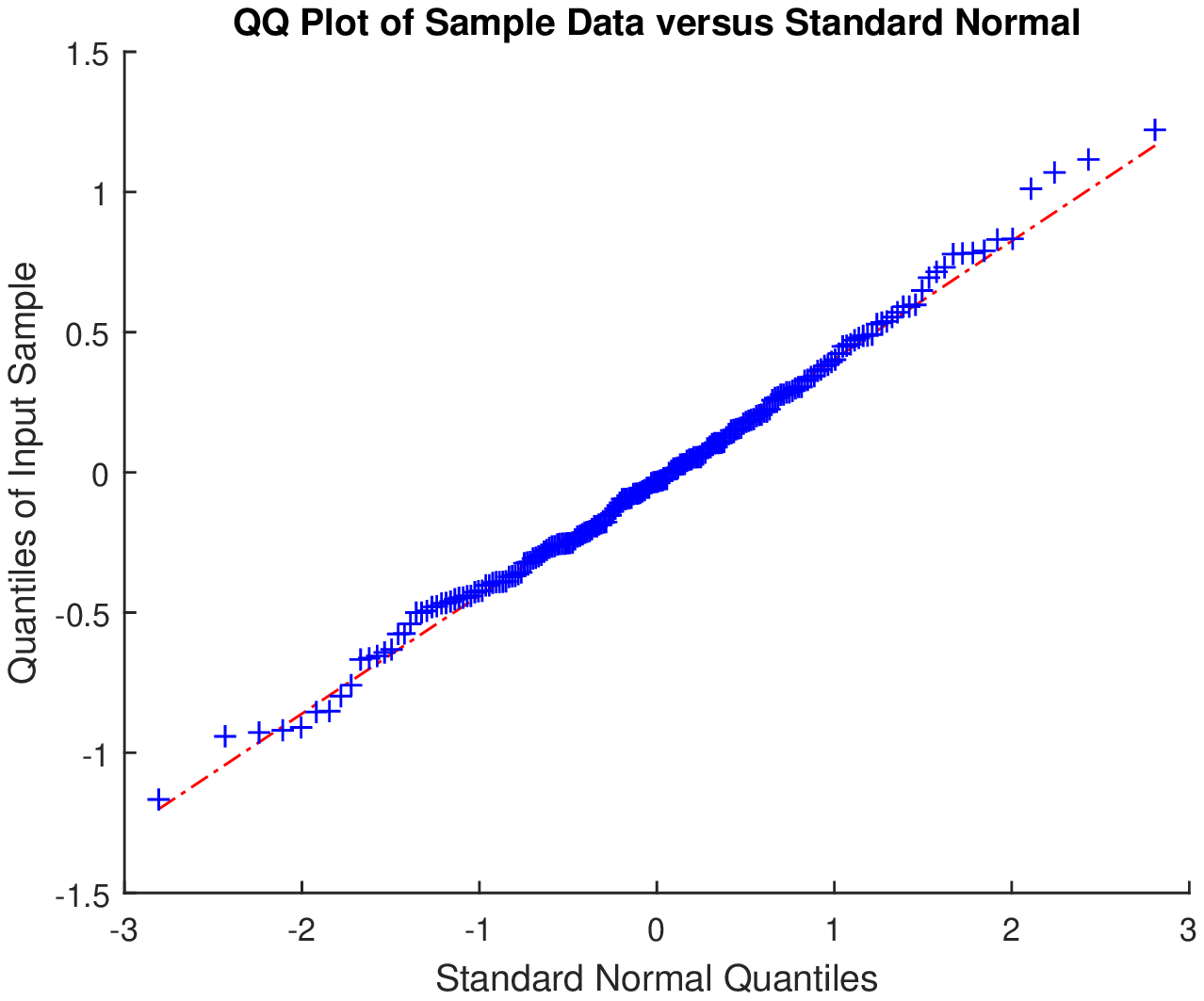}
\caption{\footnotesize{$\rho =0.9$, $TR\approx20\%$, $n= 300$, $B=200$ for $ CR\approx 0\%,$ $10\%$ and $40\%$, respectively.}}\label{Figure5.5}
 \end{center}
\end{figure}
\begin{figure}[H]
\begin{center}
 \includegraphics[scale=0.36]{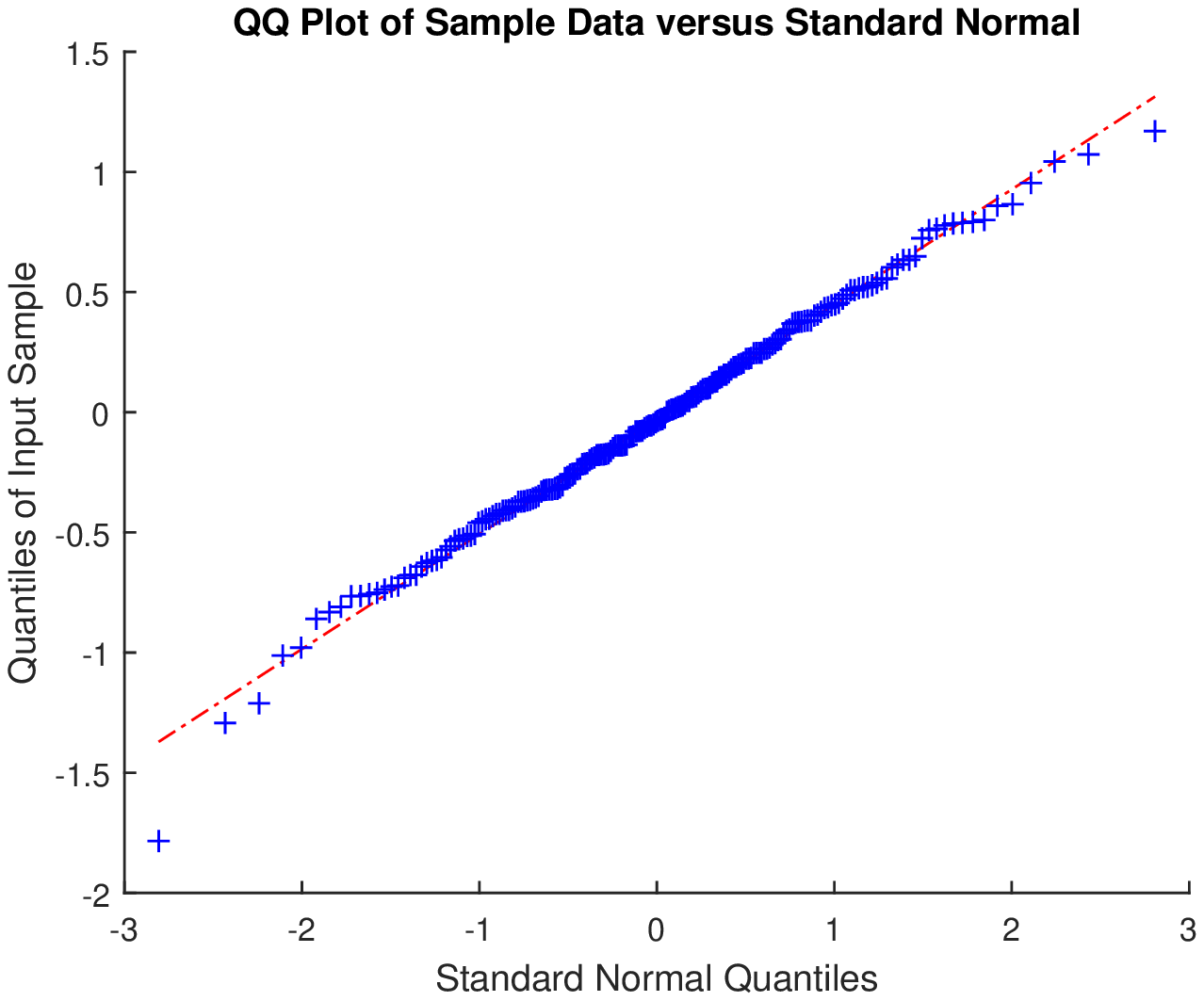}
 \includegraphics[scale=0.36]{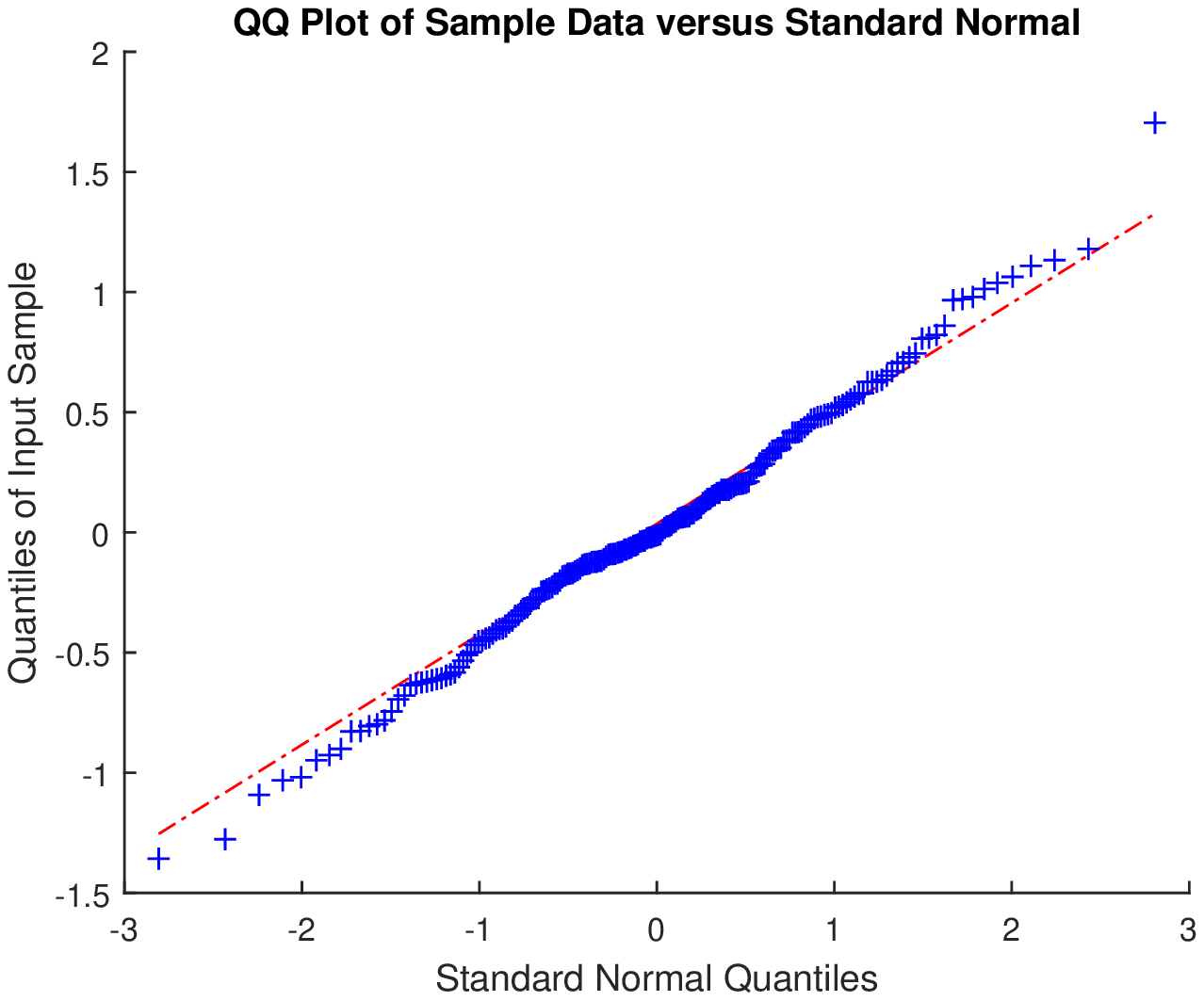}
 \includegraphics[scale=0.36]{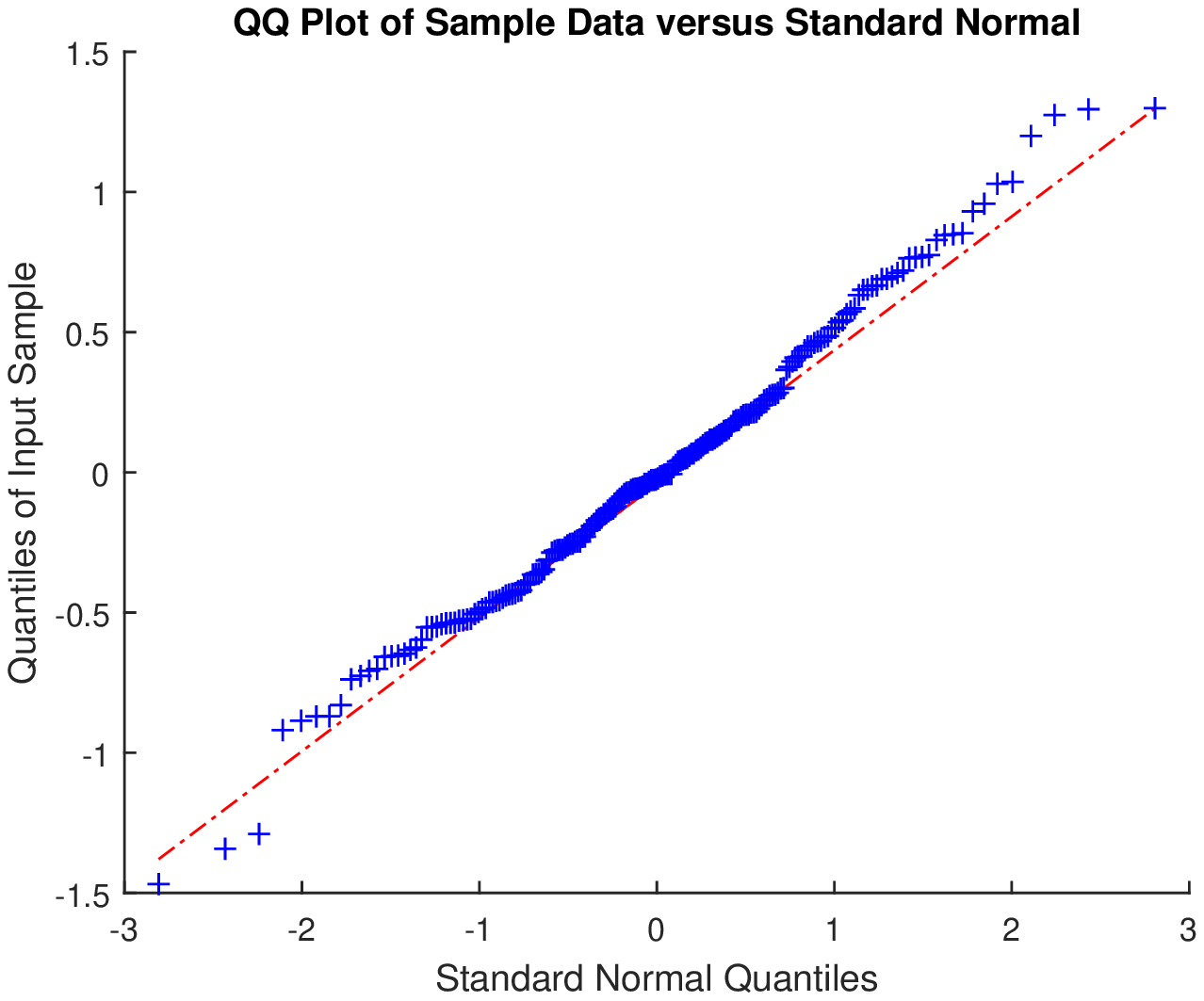}
\caption{\footnotesize{$\rho =0.9$, $CR\approx20\%$, $n= 300$, $B=200$ for $ TR\approx 0\%,$ $10\%$ and $40\%$, respectively.}}\label{Figure5.6}
 \end{center}
\end{figure}
From Figure 1-3 it is clearly seen that the Normality in the distribution of the estimator performs better for large sample sizes and seems to be not affected by high values of truncation or censoring, which is also clear from the corresponding QQ-Plot, Figure 4-6. Recall that these comparisons are given for the point $x=0$. We have also good behavior for $x\neq0$ (here we take $x\in[-1.5,1.5]$), Figure 7.
\begin{figure}[H]
\begin{center}
 \includegraphics[scale=0.6]{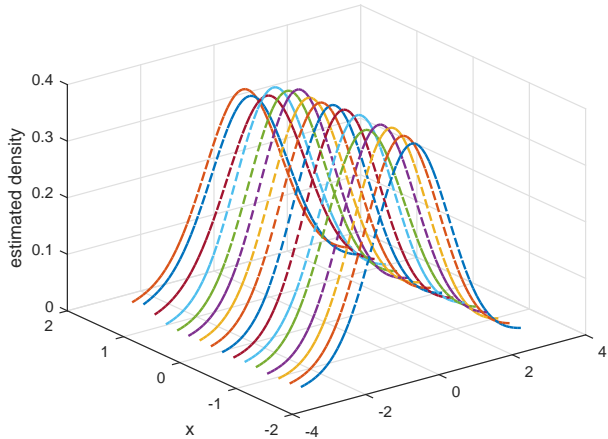}
 \caption{\footnotesize{$\rho =0.9$, $CR\approx20\%$, $ TR\approx20\%$, $n= 100$, $B=200$.}}\label{Figure5.7}
 \end{center}
\end{figure}

Next, we simulate the confidence intervals of $m(x)$ at level $95\%$. We give the plot of the median for $x\in[-2,2]$, based on $200$ replications, by choosing different sample sizes with $\rho=0.9$, $CR\approx20\%$ and $TR\approx20\%$, (Figure 8).
\begin{figure}[H]
\begin{center}
 \includegraphics[scale=0.36]{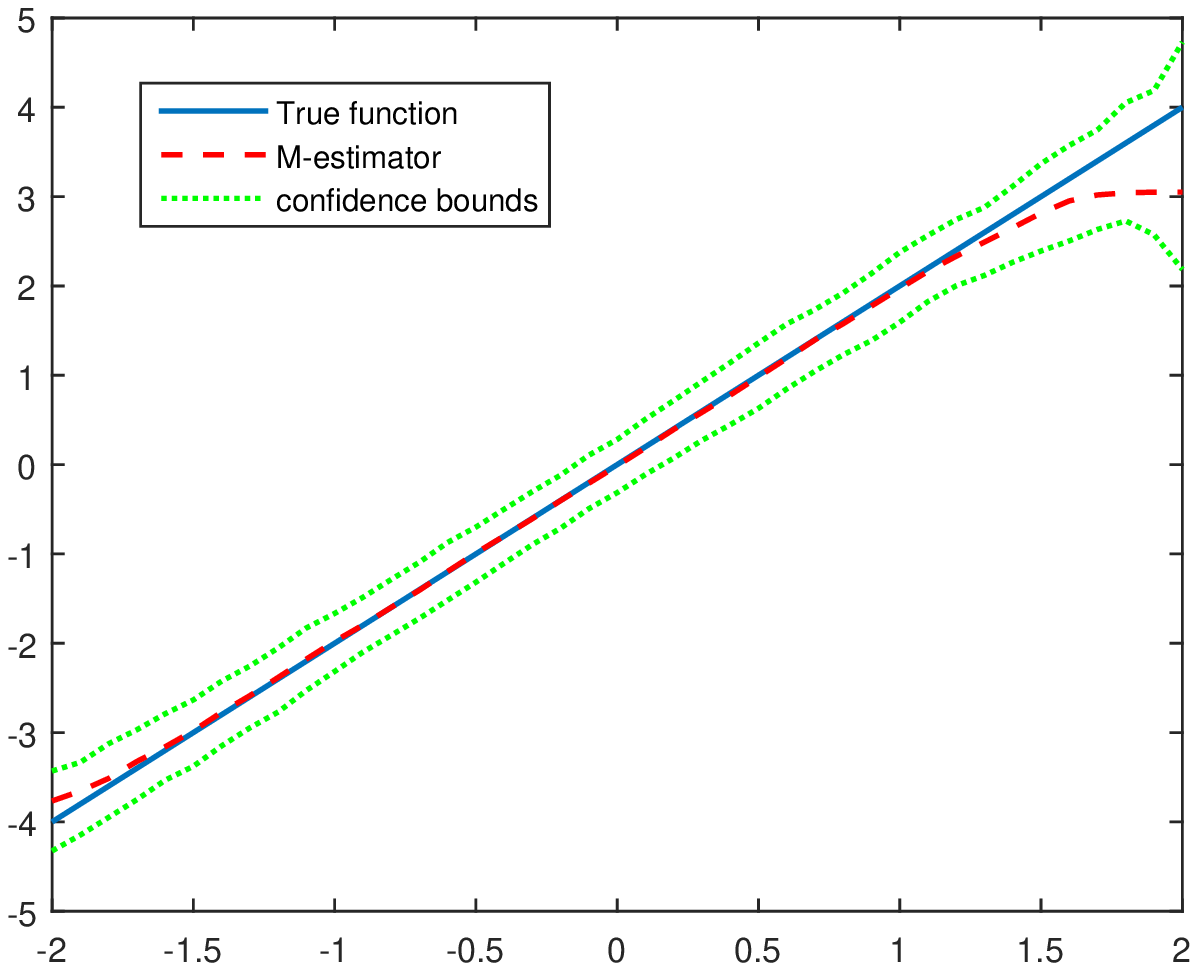}
 \includegraphics[scale=0.36]{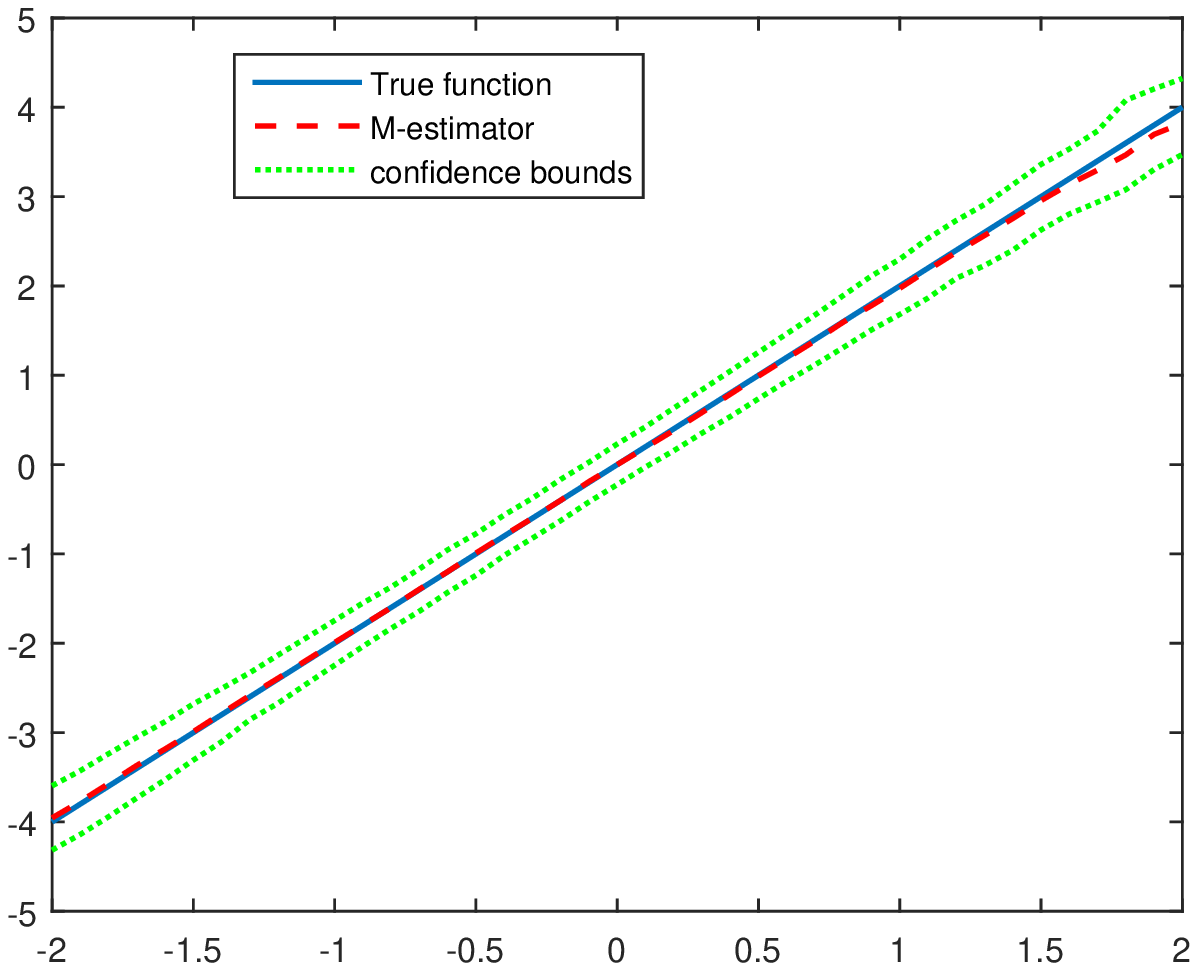}
 \includegraphics[scale=0.36]{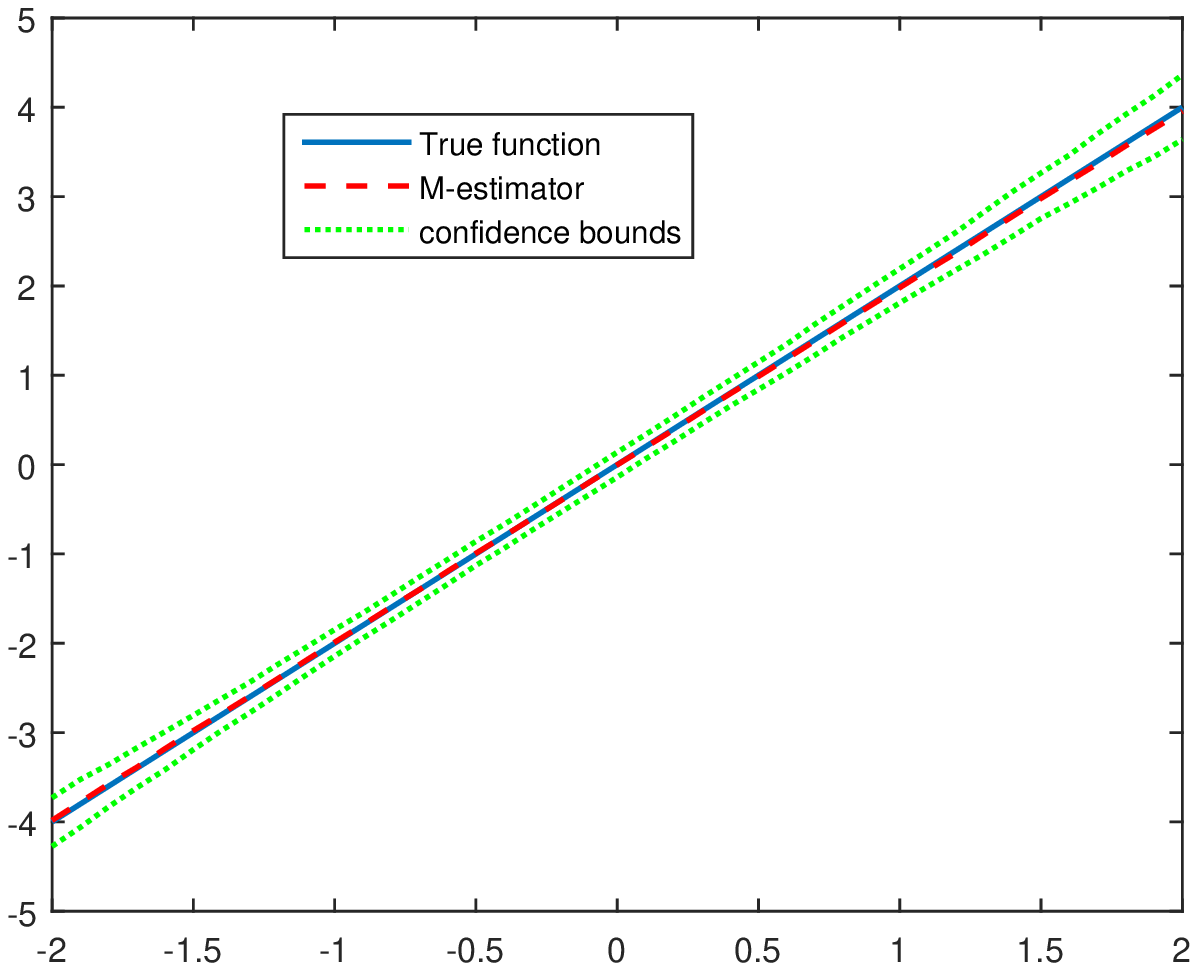}
\caption{\footnotesize{Confidence intervals of M-estimate with $\rho =0.9$, $TR\approx20\%$, $CR\approx 20\%$ for $n= 50$, $100$ and $300$ respectively.}}\label{Figure5.7}
 \end{center}
\end{figure}
Here, again we see clearly the improvement of the interval's width for a high sample sizes. Which is also confirmed in Table 1, in which we report the average widths of $95\%$ confidence intervals and the coverage probabilities based on $200$ replications. We take $x\in[-1,1]$, with increment 0.1, and different combination of $n$, $CR$ and $TR$. As we have already pointed out, the estimators quality is better when the simple sizes increases and the censoring rate decreases, while it does not seem to be affected by the percentage of truncation.\\
\begin{table}[H]
  \centering
\[
\begin{tabular}{llllllll}
\cline{3-8}
&  & $\ \ \ \ \ \ \ \ \ n=$ & $50$ & $\ \ \ \ \ \ \ \ \ n=$ & $100$ & $\ \ \
\ \ \ \ \ \ n=$ & $300$ \\ \hline
$TR\%$ & $CR\%$ & Cov Prob & A.width & Cov Prob & A.width & Cov Prob &
A.width \\ \hline
$20$ & $%
\begin{array}{c}
10 \\
40%
\end{array}%
$ & $%
\begin{array}{c}
0.7250\\
0.5450%
\end{array}%
$ & $%
\begin{array}{c}
0.7370 \\
0.9227%
\end{array}%
$ & $%
\begin{array}{c}
0.8900 \\
0.7000%
\end{array}%
$ & $%
\begin{array}{c}
0.5122 \\
0.6630%
\end{array}%
$ & $%
\begin{array}{c}
0.9620 \\
0.9250%
\end{array}%
$ & $%
\begin{array}{c}
0.2877 \\
0.3751%
\end{array}%
$ \\
$60$ & $%
\begin{array}{c}
10 \\
40%
\end{array}%
$ & $%
\begin{array}{c}
0.7642 \\
0.5213%
\end{array}%
$ & $%
\begin{array}{c}
0.7267 \\
0.9744%
\end{array}%
$ & $%
\begin{array}{c}
0.8845 \\
0.7251%
\end{array}%
$ & $%
\begin{array}{c}
0.5008 \\
0.6560%
\end{array}%
$ & $%
\begin{array}{c}
0.9515 \\
0.9213%
\end{array}%
$ & $%
\begin{array}{c}
0.2880 \\
0.3657%
\end{array}%
$ \\ \hline
\end{tabular}
\]
\caption{\footnotesize{Average widths and Coverage Probabilities of confidence intervals}\label{tab5.1}}
\end{table}
\section{Auxiliary results and proofs.}
\subsection{Weak Consistency}
The proof of Proposition \ref{Propo1} is mainly based on the following decomposition
\begin{eqnarray}\label{5.1}
\widehat{\Psi }_{x}\left( x,\theta \right) -\Psi _{x}\left( x,\theta \right)
&=&\left\{ \widehat{\Psi }_{x}\left( x,\theta \right) -\widetilde{\Psi }_{x}\left(
x,\theta \right) \right\} +\left\{ \widetilde{\Psi }_{x}\left( x,\theta \right) -%
\mathbf{E}\left( \widetilde{\Psi }_{x}\left( x,\theta \right) \right) \right\}
\notag \\
&&+\left\{ \mathbf{E}\left( \widetilde{\Psi }_{x}\left( x,\theta \right) \right)
-\Psi_{x} \left( x,\theta \right) \right\} .
\end{eqnarray}
We first give the explicit expression of the bias (Lemma \ref{lem1}) and the variance (Lemma \ref{lem2}) of $%
\widetilde{\Psi }_{x}\left( x,\theta \right) $  and next we show that $\left\{ \widehat{\Psi }_{x}\left( x,\theta \right) -\widetilde{\Psi }_{x}%
\left( x,\theta \right) \right\} $ is negligible (Lemma \ref{lem3}).
\vskip 4mm
\begin{lem}\label{lem1} Under Assumptions (K) and E1 we have
\begin{equation*}
\mathbf{E}\left( \widetilde{\Psi }\left( x,\theta \right) \right) -\Psi
\left( x,\theta \right) =\frac{h_{n}^{2}}{2}\frac{\partial ^{2}\Psi
_{x}\left( x',\theta \right) }{\partial x_{i}\partial x_{j}}%
\sum\limits_{1\leq i,j\leq d}\int\limits_{%
\mathbb{R}
^{d}}w_{i}w_{j}K\left( w\right) dw+o\left( h_{n}^{2}\right) \text{ \ \
\textit{as}}\mathit{\ }n\rightarrow \infty,
\end{equation*}
where $x'\in\mathcal{U}(x)$.
\end{lem}
{\bf Proof}  The result follows immediately from the proof of Lemma 6.1 in Benseradj and Guessoum (2020).
$\hfill\Box$

\vskip 4mm
\begin{rmk}\label{rem5.1}
Note that conditions K2 and E1 are needed to get the expression of the bias, whereas the first part of condition R1 and condition K1 suffice to get consistency.
\end{rmk}
\vskip 4mm
\begin{lem}\label{lem2} Under Assumptions K1, (R) and M we have
\begin{equation*}
Var\left( \widetilde{\Psi }_{x}\left( x,\theta \right) \right) =\frac{\mu }{%
nh_{n}^{d}}\Gamma _{x}\left( x,\theta \right) \int\limits_{%
\mathbb{R}
^{d}}K^{2}\left( w\right) dw+o(\frac{1 }{nh_{n}^{d}})\text{ \ \ \textit{as}}%
\mathit{\ }n\rightarrow \infty .
\end{equation*}
\end{lem}
{\bf Proof}  For $x\in supp(v),$ and $1\leq i\leq n,$ let
\begin{equation*}
Q_{i}\left( x,\theta \right) :=\frac{1}{\sqrt{h_{n}^{d}}}K\left( \frac{x-X_{i}}{h_{n}}\right) \psi
_{x}^{\ast }\left( Z_{i},\theta \right) -\mathbf{E}\left[\frac{1}{\sqrt{h_{n}^{d}}}K\left( \frac{%
x-X_{i}}{h_{n}}\right) \psi _{x}^{\ast }\left( Z_{i},\theta \right) \right] ,
\end{equation*}%
where $\psi _{x}^{\ast }\left( Z_{i},\theta \right):=\frac{\mu\delta_{i}\psi _{x}\left( Z_{i}-\theta \right)}{L(Z_{i})\overline{G}(Z_{i})}$.
Clearly from (\ref{Psy}) we have%
\begin{equation}\label{5.2}
\sqrt{nh_{n}^{d}}\left(\widetilde{\Psi }_{x}\left( x,\theta \right) -\mathbf{E}\left( \widetilde{%
\Psi }_{x}\left( x,\theta \right) \right)\right) =\dfrac{1}{\sqrt{n}}\sum%
\limits_{i=1}^{n}Q_{i}\left( x,\theta \right) ,\ n\geq 1.
\end{equation}%
To lighten the notations set $Q_{i}\left( x,\theta \right) =:Q_{i},$ and define the covariance term $\gamma _{j}:=Cov(Q_{1},Q_{j+1}),$ $j\geq 0.$ In view  of stationarity, we have%
\begin{eqnarray}\label{5.3}
Var\left( \sqrt{nh_{n}^{d}}\widetilde{\Psi }_{x}\left( x,\theta \right) \right)  &=&Var\left( Q_{1}\right)
+\frac{2}{n}\sum\limits_{1\leq i<j\leq n}Cov\left( Q_{i},Q_{j}\right)
\notag \\
&=:& \gamma _{0}+2\sum\limits_{j=1}^{n-1}\left( 1-\frac{j}{n}%
\right) \gamma _{j} .
\end{eqnarray}%
From (\ref{2.4}), Assumptions K1, R1, and the Dominated Convergence Theorem we have for j=0,
\begin{eqnarray}\label{5.4}
\gamma _{0} &=&\mathbf{E}\left( Q_{1}\right) ^{2}=\frac{\mu }{h_{n}^{d}}%
\int\limits_{\mathbb{R}^{d}}K^{2}\left( \frac{x-s}{h_{n}}\right) \mathbb{E}\left[ \frac{\psi
_{x}^{2}\left( Y-\theta \right) }{L\left( Y\right) \overline{G}\left(
Y\right) }|X=s\right] v\left( s\right) ds  \notag \\
&&-\frac{1}{h_{n}^{d}}\left( \int\limits_{\mathbb{R}^{d}}K\left( \frac{x-s}{h_{n}}\right) \mathbb{E}\left[ \psi _{x}\left(
Y-\theta \right) |X=s\right] v\left( s\right) ds\right) ^{2}  \notag \\
&=&\mu\int\limits_{\mathbb{R}^{d}}K^{2}\left( w\right)\mathbb{E}\left[ \frac{ \psi _{x}^{2}\left(
Y-\theta \right) }{L\left( Y\right) \overline{G}\left( Y\right) }|X=x-h_{n}w\right]
v\left( x-h_{n}w\right)dw+O\left(h_{n}^{d}\right)\notag \\
&\rightarrow&\mu\mathbb{E}\left[ \frac{ \psi _{x}^{2}\left(
Y-\theta \right) }{L\left( Y\right) \overline{G}\left( Y\right) }|X=x\right]
v\left( x\right) \int\limits_{\mathbb{R}^{d}}K^{2}\left( w\right) dw\notag \\
&=:&\mu\Gamma _{x}\left( x,\theta \right) \kappa.
\end{eqnarray}%
For $j\geq 1,$ on the one hand, using a properties of
conditional expectation and a change of variable, we get%
\begin{small}
\begin{eqnarray}\label{5.5}
\left\vert \gamma _{j}\right\vert &\leq &\frac{1}{h_{n}^{d}}\mathbf{E}\left[ K\left(\frac{ x-X_{1}}{h_{n}}\right) K\left(\frac{x-X_{j+1}}{h_{n}}\right) \left\vert \psi _{x}^{\ast }\left( Z_{1},\theta \right)\psi _{x}^{\ast }\left( Z_{j+1},\theta \right) \right\vert \right]\notag \\
&&+\frac{1}{h_{n}^{d}}\mathbf{E}^{2}\left[K\left(\frac{ x-X_{1}}{h_{n}}\right) \left\vert \psi _{x}^{\ast}\left( Z_{1},\theta \right) \right\vert \right]   \notag \\
&\leq &\frac{\mu ^{2}}{h_{n}^{d}L^{2}\left( a_{H}\right) \overline{G}^{2}\left(
b\right) }\mathbf{E}\left[  K\left(\frac{ x-X_{1}}{h_{n}}\right) K\left(\frac{x-X_{j+1}}{h_{n}}\right) \left\vert \psi _{x}\left( Y_{1}-\theta \right) \psi_{x}\left( Y_{j+1}-\theta \right) \right\vert \right]   \notag \\
&&+\frac{\mu ^{2}}{h_{n}^{d}}\mathbf{E}^{2}\left[ K\left(\frac{ x-X_{1}}{h_{n}}\right)\frac{\delta
_{1}\left\vert \psi _{x}\left( Y_{1}-\theta \right) \right\vert }{L\left(
Y_{1}\right) \overline{G}\left( Y_{1}\right) }\right]   \notag \\
&\leq &\frac{C^{2}}{h_{n}^{d}}\int\limits_{\mathbb{R}^{d}}\int\limits_{\mathbb{R}^{d}}K\left( \frac{x-w}{h_{n}}\right) K\left( \frac{x-t}{h_{n}}\right)
\mathbf{E}\left[ \left\vert \psi _{x}\left( Y_{1}-\theta \right) \psi_{x}\left( Y_{j+1}-\theta \right) \right\vert |X_{1}=w,X_{j+1}=t\right]\notag \\
&&\times \mathbf{v}_{j}\left( w,t\right) dwdt+\frac{1}{h_{n}^{d}}\left( %
\int\limits_{\mathbb{R}^{d}}K\left( \frac{x-w}{h_{n}}\right) \mathbb{E}\left[ \left\vert \psi
_{x}\left( Y-\theta \right) \right\vert |X=w\right] v\left( w\right)
dw\right) ^{2}  \notag \\
&=&O\left( h_{n}^{d}\right),
\end{eqnarray}%
\end{small}
by assumptions K1, R2 and R3. On the other hand, it follows from Davydov's Lemma (see, Ferraty and Vieu 2006) that for $s>2$, we have%
\begin{equation*}
\left\vert \gamma _{j}\right\vert \leq C\alpha \left( j\right)
^{1-2/s}\left( \mathbf{E}\left\vert Q_{j}\right\vert ^{s}\right) ^{2/s},
\end{equation*}%
where Minkowski inequality and Assumptions K1 and R2, permits to write
\begin{eqnarray*}
\mathbf{E}\left\vert Q_{j}\right\vert ^{s} &\leq &\frac{2^{s}\mu ^{s-1}}{%
h_{n}^{d\left( s/2-1\right) }L^{s-1}\left( a_{H}\right) \overline{G}%
^{s-1}\left( b\right) }\int\limits_{\mathbb{R}^{d}}K^{s}\left( w\right) \mathbb{E}\left[ \left\vert \psi _{x}\left(
Y-\theta \right) \right\vert ^{s}|X=x-h_{n}w\right] v\left( x-h_{n}w\right)
dw  \notag \\
&=&O\left( \frac{1}{h_{n}^{d\left( s/2-1\right) }}\right) ,
\end{eqnarray*}%
which yields%
\begin{equation}\label{5.7}
\left\vert \gamma _{j}\right\vert \leq c\left( h_{n}^{d}\right)
^{2/s-1}\alpha \left( j\right) ^{1-2/s}.
\end{equation}%
To end the proof of Lemma \ref{lem2} and from (\ref{5.3}) and (\ref{5.4}) it suffices to show that
\begin{equation}\label{5.8}
\sum\limits_{j=1}^{n-1}\left( 1-\dfrac{j}{n}\right) \gamma
_{j}=o\left( 1\right) ,
\end{equation}%
for this purpose let $\left\{ \zeta_{n}\right\}$ a positive sequence of integers (specified below) such that $\zeta_{n}<n$, $\zeta_{n}\rightarrow \infty $ and $\zeta_{n}h_{n}^{d}\rightarrow 0,$ and decompose the sum $\sum\limits_{j=1}^{n-1}%
\left( 1-\dfrac{j}{n}\right) \gamma _{j}$ into two terms as follows%
\begin{eqnarray}\label{5.9}
\sum\limits_{j=1}^{n-1}\left( 1-\frac{j}{n}\right) \gamma _{j}
&=&\sum\limits_{j=1}^{\zeta_{n}}\left( 1-\frac{j}{n}\right) \gamma
_{j}+\sum\limits_{j=\zeta_{n}+1}^{n-1}\left( 1-\frac{j}{n}\right) \gamma _{j}.
\end{eqnarray}%
Clearly, using (\ref{5.5}) we have
\begin{eqnarray}\label{5.10}
\sum\limits_{j=1}^{\zeta_{n}}\left( 1-\frac{j}{n}\right) \left|\gamma
_{j}\right| &\leq &\zeta_{n}\underset{1\leq
j\leq \zeta_{n}}{\max }\left\{ \left\vert \gamma _{j}\right\vert \right\}  \notag
\\
&=&O\left( h_{n}^{d}\zeta_{n}\right)  \notag \\
&=&o(1).
\end{eqnarray}%
Now for the second sum of the right hand side of (\ref{5.9}) and by (\ref{5.7}) we have
\begin{eqnarray*}
\sum\limits_{j=\zeta_{n}+1}^{n-1}\left( 1-\frac{j}{n}\right) \left|\gamma _{j}\right|
&\leq &c\left( h_{n}^{d}\right) ^{2/s-1}\sum\limits_{j=\zeta_{n}}^{\infty
}\alpha \left( j\right) ^{1-2/s} \\
&\leq &c\left( h_{n}^{d}\right) ^{2/s-1}\zeta_{n}^{-\varrho
}\sum\limits_{j=\zeta_{n}}^{\infty }j^{\varrho }\alpha \left( j\right) ^{1-2/s}
\\
&\leq &c\left( h_{n}^{d}\right) ^{2/s-1}\zeta_{n}^{-\varrho
}\sum\limits_{j=\zeta_{n}}^{\infty }j^{-\lambda \left( 1-2/s\right) +\varrho },
\end{eqnarray*}%
where $\varrho$ is such that $\left( 1-2/s\right) <\varrho <\lambda \left( 1-2/s\right)
-1$.  By choosing $\zeta_{n}=\left( h_{n}^{d}\right) ^{\left( 2/s-1\right)
/\varrho }$ we have
\begin{equation}\label{5.11}
\sum\limits_{j=\zeta_{n}+1}^{n-1}\left( 1-\frac{j}{n}\right) \left|\gamma _{j}\right|  =o(1).
\end{equation}%
Note that $\varrho $ exists thanks to Assumption M.
Hence from (\ref{5.9})-(\ref{5.11}), (\ref{5.8}) is established which completes the proof of
Lemma \ref{lem2}.
$\hfill\Box$

\vskip 4mm
\begin{lem}\label{lem3}Under Assumptions K1, R2 and M we have
\begin{equation*}
\left\vert \widehat{\Psi }\left( x,\theta \right) -\widetilde{\Psi }\left(
x,\theta \right) \right\vert =O_{\bf{P}}\left( n^{-1/2}\right)\text{ \ \ \textit{as}}%
\mathit{\ }n\rightarrow \infty  .
\end{equation*}
\end{lem}
{\bf Proof}  We have the following decomposition
\begin{eqnarray*}
\left\vert \widetilde{\Psi }\left( x,\theta \right) -\widehat{\Psi }\left(
x,\theta \right) \right\vert &\leq &\left[ \dfrac{\underset{t\in \left[
a_{H},b\right] }{\sup }\left\vert F_{n}\left( t\right) -F\left( t\right)
\right\vert }{\underset{t\in \left[ a_{H},b\right] }{\inf }\left\vert
C\left( t\right) \right\vert -\underset{t\in \left[ a_{H},b\right] }{\sup }%
\left\vert \left( C_{n}\left( t\right) -C\left( t\right) \right) \right\vert
}\right. \\
&&\left. +\dfrac{\underset{t\in \left[ a_{H},b\right] }{\sup }\left\vert
\left( C_{n}\left( t\right) -C\left( t\right) \right) \right\vert }{\underset%
{t\in \left[ a_{H},b\right] }{\inf }\left\vert C\left( t\right) \right\vert
\left( \underset{t\in \left[ a_{H},b\right] }{\inf }\left\vert C\left(
t\right) \right\vert -\underset{t\in \left[ a_{H},b\right] }{\sup }%
\left\vert \left( C_{n}\left( t\right) -C\left( t\right) \right) \right\vert
\right) }\right] \\
&&\times \left[ \frac{1}{nh_{n}^{d}}\sum\limits_{i=1}^{n}K\left( \frac{%
x-X_{i}}{h_{n}}\right) \delta _{i}\left\vert \psi _{x}\left( Z_{i}-\theta
\right) \right\vert \right] \\
&=:&I_{1n}\times I_{2n}\left( x,\theta \right).
\end{eqnarray*}
On the one hand by Lemma \ref{LemmaA5} in the Appendix we have {\small $\underset{t\in \left[ a_{H},b\right] }{\sup }\left\vert F_{n}\left( t\right) -F\left( t\right) \right\vert=O_{\bf P}(n^{-1/2})$} and {\small $\underset{t\in \left[ a_{H},b\right] }{\sup }%
\left\vert C_{n}\left( t\right) -C\left( t\right) \right\vert=O_{\bf P}(n^{-1/2})$}. Moreover, from the definition of {\small $C\left( t\right) $}, we have {\small $C\left( t\right) \geq \frac{1}{\mu}L\left( a_{H}\right) \overline{H}\left(b\right) >0,$} for all $t\in \left[ a_{H},b\right] ,$ which leads to $I_{1n}=O_{\bf P}(n^{-1/2})$ as $n\rightarrow \infty .$ In the other hand, by Markov inequality for $\epsilon >0,$ under Assumptions K1 and R2, and from (\ref{2.4}), we have
\begin{eqnarray*}
\mathbf{P}\left( I_{2n}\left( x,\theta \right) >\epsilon \right) &\leq &%
\frac{\mathbf{E}\left( I_{2n}\left( x,\theta \right) \right) }{\epsilon } \\
&\leq &\frac{1}{\mu \epsilon }\int\limits_{\mathbb{R}^{d}}K\left( w\right) \mathbb{E}\left[ \left\vert \psi _{x}\left( Y-\theta
\right) \right\vert X_{1}=x-wh_{n}\right] v(x-wh_{n})dw ,\\
&\rightarrow&\frac{1}{\mu \epsilon }\mathbb{E}\left[ \left\vert \psi _{x}\left( Y-\theta
\right) \right\vert X_{1}=x\right] v(x)
\end{eqnarray*}
which implies that $I_{2n}\left( x,\theta \right) =O_{\bf P}\left( 1\right)$ and this complete the proof of the Lemma \ref{lem3}.
$\hfill\Box$

{\bf Proof of Proposition \ref{Propo1}}  The result follow immediately from decomposition (\ref{5.1}), Tchebychev inequality, Lemmas \ref{lem1}, \ref{lem2} and \ref{lem3}.
$\hfill\Box$

{\bf Proof of Theorem \ref{theo1}}  The proof of Theorem \ref{theo1} is analogous to that of Theorem 4.1 in Benseradj and Guessoum (2020), therefore, it is omitted.
$\hfill\Box$
\subsection{Asymptotic Normality}
{\bf Proof of Proposition \ref{Propo2}} Lemma \ref{lem1} and Lemma \ref{lem3} under Assumption B1 give
\begin{equation*}
\sqrt{nh_{n}^{d}}\left( \left\{ \widehat{\Psi }_{x}\left( x,\theta \right) -%
\widetilde{\Psi }_{x}\left( x,\theta \right) \right\} +\left\{ \mathbf{E}%
\left( \widetilde{\Psi }_{x}\left( x,\theta \right) \right) -\Psi _{x}\left(
x,\theta \right) \right\} \right) =o_{\bf P}\left( 1\right) ,
\end{equation*}%
hence to prove Proposition \ref{Propo2} and by (\ref{5.1}) it suffices to prove that $\dfrac{%
\sqrt{nh_{n}^{d}}\left[ \widetilde{\Psi }\left( x,\theta \right) -\mathbf{E}%
\left( \widetilde{\Psi }\left( x,\theta \right) \right) \right] }{\sigma
_{0}\left( x,\theta \right) }$ converge in distribution to standard normal, which is given by the following Lemma.
$\hfill\Box$

\vskip 4mm
\begin{lem}\label{lem4}Under Assumptions K1, (R), M and B2 we have
\begin{equation*}
\dfrac{\sqrt{nh_{n}^{d}}\left[ \widetilde{\Psi }\left( x,\theta \right) -%
\mathbf{E}\left( \widetilde{\Psi }\left( x,\theta \right) \right) \right] }{%
\sigma _{0}\left( x,\theta \right) }\overset{D}{\rightarrow }\mathcal{N}\left(
0,1\right) \text{ \ \ \textit{as}}%
\mathit{\ }n\rightarrow \infty  .
\end{equation*}
\end{lem}
{\bf Proof}  We employ the Doob's small-block and large-block technique (see Doob 1953). Namely, divide the set $\left\{ 1,2,...,n\right\} $ into $2w_{n}+1$ subsets with large-block of size $p:=p_{n}$ and small block of size $q:=q_{n},$ where $%
p_{n},$ $q_{n}$ are integers tending to infinity along with $n.$ Set
\begin{equation*}
w:=w_{n}=\left[ \frac{n}{p+q}\right] ,
\end{equation*}
where $\left[ x\right] $ stands for the integers part of $x$. Note that condition B2 implies that there exists a sequence of positive
integers $\left\{ \xi _{n}\right\} ,$ $\xi _{n}\rightarrow \infty $ along with $n$, such that
\begin{equation}\label{5.24}
\xi _{n}q_{n}=o\left( \sqrt{nh_{n}^{d}}\right) \text{ \ } \text{and} \text{ \ } \xi _{n}\left(n/h_{n}^{d}\right)^{1/2}\alpha\left(q_{n}\right)=o(1) \text{ \ }as \text{ \ }n\rightarrow \infty.
\end{equation}
Now define the small block size as $q_{n}=\left[\frac{\sqrt{{nh_{n}^{d}}}}{\xi _{n}logn}\right]$ and the large block size as $p_{n}=\left[\frac{\sqrt{{nh_{n}^{d}}}}{\xi _{n}\sqrt{logn}}\right].$ Then (\ref{5.24}) and simple algebra show that
\begin{equation}\label{5.18}
q/p\longrightarrow 0,\text{ \ }w\alpha \left( q\right) \longrightarrow 0,%
\text{ \ }wq/n\longrightarrow 0,\ \ p/n\longrightarrow 0,\text{ \ }%
pw/n\longrightarrow 1, \text{ \ } \text{as}\text{ \ }n\rightarrow \infty.
\end{equation}
For $1\leq z\leq w$, define the random variables
\begin{equation}\label{5.13}
\chi _{zn}\left( x\right) =\sum\limits_{i=k_{z}}^{k_{z}+p}Q_{i}\left(
x, \theta\right) ,\text{ \ }\chi _{zn}^{\prime }\left( x\right)
=\sum\limits_{i=l_{z}}^{l_{z}+q-1}Q_{i}\left( x,\theta\right) ,\text{ \ }\chi
_{wn}^{\prime \prime }\left( x\right) =\sum\limits_{i=w\left(p+q\right)
+1}^{n}Q_{i}\left( x,\theta\right) ,
\end{equation}
where $k_{z}=\left( z-1\right) \left( p+q\right) +1;$ $\ \ l_{z}=\left(
z-1\right) \left( p+q\right) +p+1;$ then from (\ref{5.2}), we have
\begin{eqnarray*}
\sqrt{nh_{n}^{d}}\left[ \widetilde{\Psi }\left( x,\theta \right) -%
\mathbf{E}\left( \widetilde{\Psi }\left( x,\theta \right) \right) \right]
&=&n^{-1/2}\left\{\sum\limits_{z=1}^{w}\sum\limits_{i=k_{z}}^{k_{z}+p}Q _{i}\left(
x,\theta\right) +\sum\limits_{z=1}^{w}\sum\limits_{i=l_{z}}^{l_{z}+q-1}Q
_{i}\left( x,\theta\right) +\sum\limits_{i=w\left( p+q\right) +1}^{n}Q
_{i}\left( x,\theta\right) \right\} \\
&=&n^{-1/2}\left\{ \sum\limits_{z=1}^{w}\chi _{zn}\left( x\right)
+\sum\limits_{z=1}^{w}\chi _{zn}^{\prime }\left( x\right) +\chi
_{wn}^{\prime \prime }\left( x\right) \right\} \\
&=:&n^{-1/2}\left\{ s_{n}^{\prime }\left( x\right) +s_{n}^{\prime \prime
}\left( x\right) +s_{n}^{\prime \prime \prime }\left( x\right) \right\} .
\end{eqnarray*}%
Thus we need to prove that as $n\rightarrow \infty $,%
\begin{equation}\label{5.14}
n^{-1}\mathbf{E}\left( s_{n}^{\prime \prime 2}\left( x\right) \right)
\longrightarrow 0\ ;\ \ \ n^{-1}\mathbf{E}\left( s_{n}^{\prime \prime \prime
2}\left( x\right) \right) \longrightarrow 0,
\end{equation}%
\begin{equation}\label{5.15}
n^{-1}\sum\limits_{z=1}^{w}\mathbf{E}\left( \chi _{zn}^{2}\left(
x\right) \right)\longrightarrow \sigma _{0}^{2}\left( x,\theta
\right),
\end{equation}%
\begin{equation}\label{5.16}
\left\vert \mathbf{E}\left( \exp \left( itn^{-1/2}s_{n}^{\prime }\left( x\right)
\right) \right) -\prod\limits_{z=1}^{w}\mathbf{E}\left( \exp \left( itn^{-1/2}\chi
_{zn}\left( x\right) \right) \right) \right\vert \longrightarrow 0,
\end{equation}%
\begin{equation}\label{5.17}
A_{n}\left( \varepsilon \right) =\frac{1}{n}\sum\limits_{z=1}^{w}\mathbf{E}%
\left[ \chi _{zn}^{2}\left( x\right) \mathbb{I}_{\left\{ \left\vert \chi
_{zn}\left( x\right) \right\vert >\varepsilon\sigma _{0}\left( x,\theta
\right) \sqrt{n}\right\} }\right]
\longrightarrow 0\text{ \ \ \ \ \ }\forall \varepsilon >0.
\end{equation}
Indeed, the result (\ref{5.14}) implies that $s_{n}^{\prime \prime }\left( x\right) $ and $s_{n}^{\prime \prime \prime }\left( x\right) $ are asymptotically negligible in probability; the result (\ref{5.16}) shows that the summands $\chi _{zn}\left( x\right) $
in $s_{n}^{\prime }\left( x\right) $ are asymptotically independent; while the results (\ref{5.15}) and (\ref{5.17}) are the standard Lindeberg-Feller conditions for asymptotic normality of $s_{n}^{\prime }\left( x\right) $ for independent setup.\\
\textbf{1) }Let us first establish (\ref{5.14}). Observe that
\begin{eqnarray}\label{5.19}
n^{-1}\textbf{E}\left( s_{n}^{\prime \prime 2}\left( x\right) \right)
&=&n^{-1}Var\left( s_{n}^{^{\prime \prime }}\left( x\right) \right)  \notag
\\
&=&\frac{1}{n}Var\left( \sum\limits_{z=1}^{w}\chi _{zn}^{\prime }\left(
x\right) \right)  \notag \\
&=&\frac{1}{n}\sum\limits_{z=1}^{w}Var\left( \chi _{zn}^{\prime }\left(
x\right) \right) +\frac{2}{n}\sum\limits_{1\leq z_{1}<z_{2}\leq w}Cov\left(
\chi _{z_{1}n}^{\prime }\left( x\right) ,\chi _{z_{2}n}^{\prime }\left(
x\right) \right)  \notag \\
&=:&J_{1n}(x)+J_{2n}(x).
\end{eqnarray}
It follows from stationarity and Lemma \ref{lem2} that%
\begin{eqnarray}
J_{1n}(x) &=&\frac{w}{n}Var\left( \chi _{1n}^{\prime }\left( x\right)
\right)\notag \\
&=&\frac{w}{n}Var\left( \sum\limits_{i=1}^{q}Q_{i}\left( x,\theta \right) \right)\notag \\
&=&\frac{wq}{n}Var\left(\frac{1}{\sqrt{q}}\sum\limits_{i=1}^{q}Q_{i}\left( x,\theta \right) \right).\notag
\end{eqnarray}
Note that Lemma \ref{lem2} remains true when substituting $n$ by $q_{n},$ so by condition (\ref{5.18}) we get
\begin{eqnarray}\label{5.20}
J_{1n}(x) &=&\frac{wq}{n}\sigma _{0}^{2}\left( x,\theta \right)\left( 1+o\left( 1\right) \right)  \notag \\
&=&o(1).
\end{eqnarray}
For $J_{2n}(x),$ we have
\begin{eqnarray}\label{5.21}
\left\vert J_{2n}(x)\right\vert &\leq &2n^{-1}\sum\limits_{1\leq i<j\leq
n}\left\vert Cov\left( Q _{i}\left( x, \theta\right) ,Q _{j}\left( x, \theta\right)
\right) \right\vert  \notag \\
&=&o\left( 1\right),
\end{eqnarray}
by (\ref{5.8}). Hence, by (\ref{5.19})-(\ref{5.21}) we have
\begin{equation}\label{5.22}
n^{-1}\textbf{E}\left( s_{n}^{\prime \prime 2}\left( x\right) \right) =o\left(
1\right) .
\end{equation}
Similarly, it follows from stationarity, (\ref{5.4}), (\ref{5.21}) and condition (\ref{5.18})
that%
\begin{eqnarray}\label{5.23}
n^{-1}\textbf{E}\left( s_{n}^{\prime \prime \prime 2}\left( x\right) \right)  &=&%
\frac{1}{n}Var\left( \chi _{wn}^{\prime \prime }\left( x\right) \right) \notag \\
&=&\frac{1}{n}Var\left( \sum\limits_{i=w\left( p+q\right) +1}^{n}Q_{i}\left( x,\theta \right)\right) \notag \\
&=&\frac{ n-w\left( p+q\right)}{n}Var\left(Q_{1}\left( x,\theta \right)
\right)+\frac{2}{n}\sum\limits_{w(p+q)+1\leq i<j\leq n}Cov(Q_{i},Q_{j})\notag \\
&\leq &\frac{ n-w\left( p+q\right)}{n}Var\left(Q_{1}\left( x,\theta \right)\right)
+\frac{2}{n}\sum\limits_{1\leq i<j\leq n}\left\vert Cov\left(Q_{i},Q_{j}\right)
\right\vert \notag \\
&=& o\left(1\right)\left(\sigma_{0}^{2}\left(x,\theta \right)+o\left(1\right)\right)+o(1)\rightarrow 0,\text{ \ } as \text{ \ } n \rightarrow\infty,
\end{eqnarray}
therefore by (\ref{5.22}) and (\ref{5.23}), (\ref{5.14}) is proved.\\
\textbf{2)} As for (\ref{5.14}), by stationarity, condition (\ref{5.18}) and Lemma \ref{lem2}, it is easily seen that
\begin{eqnarray*}
n^{-1}\sum\limits_{z=1}^{w}\mathbf{E}\left( \chi _{zn}^{2}\left( x\right)
\right) &=&\frac{w}{n}\mathbf{E}\left( \chi _{1n}^{2}\left( x\right) \right)
\\
&=&\frac{w}{n}Var\left( \sum\limits_{i=1}^{p}Q_{i}\left( x,\theta \right) \right) \\
&=&\frac{wp}{n}Var\left( \frac{1}{\sqrt{p}}\sum\limits_{i=1}^{p}Q_{i}\left( x,\theta \right)
\right) \\
&=&\frac{p w}{n}\sigma _{0}^{2}\left( x,\theta \right)\left( 1+o\left( 1\right) \right) \rightarrow \sigma _{0}^{2}\left( x,\theta \right),
\end{eqnarray*}
and (\ref{5.15}) is proved.\\
\textbf{3)} In order to establish (\ref{5.16}), we make use of Lemma 1.1 in Volkonskii and Rozanov (1959) applied to $V_{z}=:\exp \left( itn^{-1/2}\chi_{zn}\left( x\right) \right)$ to obtain\\
\begin{equation*}
\left\vert \mathbf{E}\left( \exp \left( itn^{-1/2}s_{n}^{\prime }\left( x\right)
\right) \right) -\prod\limits_{z=1}^{w}\mathbf{E}\left( \exp \left( itn^{-1/2}\chi
_{zn}\left( x\right) \right) \right) \right\vert \leq 16\left( w-1\right)
\alpha \left( q+1\right) \leq 16w\alpha \left( q\right) ,
\end{equation*}
which tends to zero as $n$ goes to infinity thanks to condition (\ref{5.18}).\\
\textbf{4)} Finally, we establish (\ref{5.17}). To this end we need to employ a truncation method since $\psi _{x}$ is not necessarily bounded. Let
\begin{equation*}
\psi _{\xi _{n}}^{\ast }\left( Z_{i},\theta \right) :=\psi _{x}^{\ast
}\left( Z_{i},\theta \right) \mathbb{I}_{\left\{ \left\vert \psi _{x}\left(
Z_{i}-\theta \right) \right\vert \leq\xi _{n}\right\} },
\end{equation*}
and define
\begin{eqnarray}\label{5.25}
 Q_{\xi _{n},i}\left( x,\theta \right) &:=&\frac{1}{\sqrt{h_{n}^{d}}}K\left(\frac{ x-X_{i}}{h_{n}}\right)\psi_{\xi _{n}}^{\ast }\left( Z_{i},\theta \right)- \mathbf{E}\left[\frac{1}{\sqrt{h_{n}^{d}}}K\left(\frac{ x-X_{i}}{h_{n}}\right)\psi _{\xi _{n}}^{\ast }\left( Z_{i},\theta\right) \right] ,\notag\\
\widetilde{\Psi }_{x,\xi _{n}}\left( x,\theta\right) &:=&\frac{1}{nh_{n}^{d}}\text{ }\sum\limits_{i=1}^{n}K\left( \frac{%
x-X_{i}}{h_{n}}\right) \psi _{\xi _{n}}^{\ast }\left( Z_{i},\theta \right) .
\end{eqnarray}
Clearly we have
\begin{equation*}
\sqrt{nh_{n}^{d}}\left[ \widetilde{\Psi }_{x,\xi _{n}}\left( x,\theta
\right) -\mathbf{E}\left( \widetilde{\Psi }_{x,\xi _{n}}\left( x,\theta
\right) \right) \right] =n^{-1/2}\text{ }\sum\limits_{i=1}^{n}Q _{\xi
_{n},i}\left( x, \theta\right) .
\end{equation*}
Let $\chi _{z,\xi _{n}}$ be given by (\ref{5.13}) with $Q _{i}$
replaced by $Q _{\xi _{n},i}.$ It is clearly seen from (\ref{5.25}) and condition K1 that $Q _{\xi _{n},i}$ is
bounded since $\left\vert Q_{\xi _{n},i}\right\vert \leq c\frac{\xi _{n}}{%
\sqrt{h_{n}^{d}}}.$ Thus by (\ref{5.13}) we have
\begin{equation*}
\underset{1\leq z\leq w}{\max }\left\vert \chi _{z,\xi _{n}}\right\vert /%
\sqrt{n}\leq c\frac{p_{n}\xi _{n}}{\sqrt{nh_{n}^{d}}}\rightarrow 0,
\end{equation*}
by the definition of $p_{n}$. Hence for $n$ large enough, the set $\left\{ \left\vert
\chi _{z,\xi _{n}}\right\vert >\varepsilon \sigma _{0,\xi _{n}}\left(
x,\theta \right) \sqrt{n}\right\} $ becomes an empty set, where we have as in
Lemma \ref{lem2}
\begin{equation*}
Var\left(\sqrt{nh_{n}^{d}} \widetilde{\Psi }_{x,\xi _{n}}\left( x,\theta
\right) \right) \rightarrow \mu \Gamma _{x,\xi _{n}}\left( x,\theta \right)
\int\limits_{\mathbb{R}^{d}}K^{2}\left( w\right) dw=:\sigma _{0,\xi _{n}}^{2}\left(
x,\theta \right) ,
\end{equation*}
\begin{equation}\label{gamma}
\Gamma _{x,\xi _{n}}\left( x,\theta \right) :=\mathbb{E}\left[ \dfrac{
\psi _{x}^{2}\left( Y_{1}-\theta \right) \mathbb{I}_{\left\{ \left\vert \psi
\left( Y_{1}-\theta \right) \right\vert \leq \xi _{n}\right\} }}{L\left(
Y_{1}\right) \overline{G}\left( Y_{1}\right) }|X_{1}=x\right] v\left(
x\right),
\end{equation}
(compare with $\Gamma _{x}\left( x,\theta \right) $ defined in condition R1) and then (\ref{5.17}) holds for $\chi _{z,\xi _{n}}.$ Consequently we have
\begin{equation*}
n^{-1/2}\sum\limits_{i=1}^{n}Q_{\xi _{n},i}\left( x\right)
\overset{D}{\rightarrow }\mathcal{N}\left( 0,\sigma _{0,\xi _{n}}^{2}\left( x,\theta
\right) \right) .
\end{equation*}
In order to generalize the proof to $n^{-1/2}$ $\sum\limits_{i=1}^{n}Q_{i}$, following the proof of (5.35d) in Masry (2005), it suffices to show that
\begin{equation}\label{5.26}
\frac{1}{n}Var\sum\limits_{i=1}^{n}\overline{Q}_{\xi _{n},i}\left(
x\right) \rightarrow 0\text{ as }n\rightarrow \infty \text{ and then }\xi
_{n}\rightarrow \infty,
\end{equation}
where $\overline{Q}_{\xi _{n},i}:=Q _{i}-Q _{\xi _{n},i}.$\\
Note that $\overline{Q }_{\xi _{n},i}$ has the same structure as $Q_{i}$ with $\psi _{x}\left( Z_{i}-\theta \right) \mathbb{I}_{\left\{ \left\vert \psi_{x}\left( Z_{i}-\theta \right) \right\vert >\xi _{n}\right\} }$ in place of $\psi _{x}\left( Z_{i}-\theta \right).$ Hence by
the argument of Lemma \ref{lem2} we have
\begin{equation*}
\frac{1}{n}Var\sum\limits_{i=1}^{n}\overline{Q}_{\xi _{n},i}\left(
x\right) \rightarrow \mu \overline{\Gamma }_{x,\xi _{n}}\left( x,\theta \right)
\int\limits_{\mathbb{R}^{d}}K^{2}\left( w\right) dw,
\end{equation*}
where
\begin{equation*}
\overline{\Gamma }_{x,\xi _{n}}\left( x,\theta \right) :=\mathbb{E}\left[
\dfrac{ \psi _{x}^{2}\left( Y_{1}-\theta \right) \mathbb{I}_{\left\{
\left\vert \psi \left( Z_{i},\theta \right) \right\vert > \xi
_{n}\right\} }}{L\left( Y_{1}\right) \overline{G}\left( Y_{1}\right) }%
|X_{1}=x\right] v\left( x\right).
\end{equation*}
By the dominated convergence theorem $\overline{\Gamma }_{x,\xi _{n}}\left(x,\theta \right) \rightarrow 0$ as $\xi _{n}\rightarrow \infty .$ This
establishes (\ref{5.26}) and completes the proof of Lemma \ref{lem4} and the Proposition \ref{Propo2}.
$\hfill\Box$

{\bf Proof of Theorem \ref{theo2}}  Using formula (\ref{3.1}) one gets%
\begin{equation}\label{5.27}
\sqrt{nh_{n}^{d}}\left( \widehat{m}\left( x\right) -m\left( x\right) \right)
=-\dfrac{\sqrt{nh_{n}^{d}}\widehat{\Psi }_{x}\left( x,m\left( x\right)
\right) }{\dfrac{\partial \widehat{\Psi }}{\partial \theta }\left( x,%
\widehat{m}^{\ast }\left( x\right) \right) },
\end{equation}
if the denominator does not vanish. By Slutsky's Theorem and Proposition \ref{Propo2} with $\theta =m\left( x\right) $ it suffices to show that the denominator in (\ref{5.27}) converges in probability to $\dfrac{\partial \Psi }{\partial \theta }\left( x,m\left( x\right) \right) $,
which is given by the following Lemma.
$\hfill\Box$

\vskip 4mm
\begin{lem} \label{lem5} Under the assumptions of Proposition \ref{Propo1} and Conditions (R$^{\prime }$) we have
\begin{equation*}
\dfrac{\partial \widehat{\Psi }_{x}}{\partial \theta }\left( x,\widehat{m}%
^{\ast }\left( x\right) \right) \overset{\bf{P}}{\rightarrow }\dfrac{\partial
\Psi _{x}}{\partial \theta }\left( x,m\left( x\right) \right) \text{ \ as }%
n\rightarrow \infty .
\end{equation*}
\end{lem}
{\bf Proof}  We write
\begin{eqnarray*}
\left\vert \dfrac{\partial \widehat{\Psi }_{x}}{\partial \theta }\left( x,%
\widehat{m}^{\ast }\left( x\right) \right) -\dfrac{\partial \Psi _{x}}{%
\partial \theta }\left( x,m\left( x\right) \right) \right\vert  &\leq
&\left\vert \dfrac{\partial \widehat{\Psi }_{x}}{\partial \theta }\left( x,%
\widehat{m}^{\ast }\left( x\right) \right) -\dfrac{\partial \widehat{\Psi }%
_{x}}{\partial \theta }\left( x,m\left( x\right) \right) \right\vert  \\
&&+\left\vert \dfrac{\partial \widehat{\Psi }_{x}}{\partial \theta }\left(
x,m\left( x\right) \right) -\dfrac{\partial \Psi _{x}}{\partial \theta }%
\left( x,m\left( x\right) \right) \right\vert  \\
&\leq &J_{1}(x)+J_{2}(x).
\end{eqnarray*}
For $J_{1}\left( x\right) $ we write
\begin{eqnarray*}
J_{1}(x) &\leq &\frac{\mu _{n}}{nh_{n}^{d}}\sum\limits_{i=1}^{n}K\left(
\frac{x-X_{i}}{h_{n}}\right) \frac{\delta _{i}}{L_{n}\left( Z_{i}\right)
\left( 1-G_{n}\left( Z_{i}\right) \right) }\left\vert \dfrac{\partial \psi
_{x}\left( Z_{i}-\widehat{m}^{\ast }\left( x\right) \right) }{\partial
\theta }-\dfrac{\partial \psi _{x}\left( Z_{i}-m\left( x\right) \right) }{%
\partial \theta }\right\vert \\
&\leq &\mu _{n}\left( L\left( a_{H}\right) -\underset{t\in \left[ a_{H},b%
\right] }{\sup }\left\vert \left( L_{n}\left( t\right) -L\left( t\right)
\right) \right\vert \right) ^{-1}\left( \overline{G}\left( b\right) -%
\underset{t\in \left[ a_{H},b\right] }{\sup }\left\vert \left( G_{n}\left(
t\right) -G\left( t\right) \right) \right\vert \right) ^{-1} \\
&&\times \underset{t\in \left[ a_{H},b\right] }{\sup }\left\vert \dfrac{%
\partial \psi _{x}\left( t-\widehat{m}^{\ast }\left( x\right) \right) }{%
\partial \theta }-\dfrac{\partial \psi _{x}\left( t-m\left( x\right) \right)
}{\partial \theta }\right\vert \frac{1}{nh_{n}^{d}}\sum\limits_{i=1}^{n}K%
\left( \frac{x-X_{i}}{h_{n}}\right) ,
\end{eqnarray*}
by Markov inequality, Assumptions K1 and R2 we can easily show that
\begin{equation*}
\frac{1}{nh_{n}^{d}}\sum\limits_{i=1}^{n}K\left( \frac{x-X_{i}}{h_{n}}%
\right) =O_{\bf{P}}\left( 1\right) .
\end{equation*}
Moreover by the uniform continuity of  $\partial/\partial \theta \left(\psi _{x}\left(
t-\theta \right)\right)$ at  $m\left( x\right) ,$ the make use
of Theorem \ref{theo1} ensures that {\small $\underset{t\in \left[ a_{H},b\right] }{\sup }%
\left\vert \dfrac{\partial \psi _{x}\left( t-\widehat{m}^{\ast }\left(
x\right) \right) }{\partial \theta }-\dfrac{\partial \psi _{x}\left(
t-m\left( x\right) \right) }{\partial \theta }\right\vert $} is negligible in
probability. Hence, Lemma 5.5 of Liang (2011) and Lemma \ref{LemmaA5} in Appendix allow to conclude that $J_{1}\left( x\right) =o_{\bf{P}}\left( 1\right)$. Finally by using similar argument as in the proof of Proposition \ref{Propo1}, under assumptions R$^{\prime }$1, R$^{\prime }$2 and R$^{\prime }$3 in place of R1, R2, and R3, we can show that $J_{2}(x)$ converges in probability to zero. Which completes the proof of Lemma \ref{lem5} and Theorem \ref{theo2}.
$\hfill\Box$
\section*{\textbf{Appendix}}
\appendix
In this section we give some supplementary new results which helped us to get our results. We will obtain a uniform weak convergence rate for the TJW product-limit estimator of the lifetime and censored distribution under dependence, which are useful tools for our study and other LTRC strong mixing framework.
\vskip 4mm
\begin{lemma}\label{LemmaA5} Suppose that $\alpha \left( k\right) =O\left(
k^{-\lambda }\right) ,$ for some $\lambda >3.$ Then for any $b\in \lbrack a_{H},b_{H})$ we have
\begin{description}
\item[i.] $\underset{y\in \left[ a_{H},b\right] }{\sup }\left\vert
C_{n}\left( y\right) -C\left( y\right) \right\vert =O_{\bf{P}}\left(
n^{-1/2}\right) ;$
\item[ii.]$\underset{y\in \left[ a_{H},b\right] }{\sup }\left\vert
F_{n}\left( y\right) -F\left( y\right) \right\vert =O_{\bf{P}}\left(
n^{-1/2}\right) ;$
\item[iii.]  $\underset{y\in \left[
a_{H},b\right] }{\sup }\left\vert G_{n}\left( y\right) -G\left( y\right)
\right\vert =O_{\bf{P}}\left( n^{-1/2}\right) ;$

\item[iv.] $\left\vert \mu _{n}-\mu \right\vert =O_{\bf{P}}\left( n^{-1/2}\right) .$
\end{description}
\end{lemma}
{\bf Proof of Lemma \ref{LemmaA5}}  To prove (i), observe that
\begin{equation}\label{6.1}
C_{n}\left( y\right) -C\left( y\right) =\left( \mathbf{L}_{n}\left( y\right)
-\mathbf{L}\left( y\right) \right) -\left( \mathbf{H}_{n}\left( y\right) -%
\mathbf{H}\left( y\right) \right) .
\end{equation}
From (6.5) in Liang (2011) we have
\begin{equation}\label{6.2b}
\underset{y\in \left[ a_{H},b\right] }{\sup }\left\vert \mathbf{L}_{n}\left(
y\right) -\mathbf{L}\left( y\right) \right\vert =O_{\bf{P}}\left( n^{-1/2}\right)
\end{equation}
By using similar arguments as in (\ref{6.2b}) we can show that
\begin{equation}\label{6.2}
\underset{y\in \left[ a_{H},b\right] }{\sup }\left\vert \mathbf{H}_{n}\left(
y\right) -\mathbf{H}\left( y\right) \right\vert =O_{\bf{P}}\left( n^{-1/2}\right)
\end{equation}
then (i) is established from (\ref{6.2b}) and (\ref{6.2}) combined with (\ref{6.1}).\\
(ii) As it is known (see, e.g. Guessoum and Tatachak (2020)), the cumulative
hazard function $\Lambda $ is defined by
\begin{equation*}
\Lambda \left( y\right) =\int\limits_{a_{H}}^{y}\frac{d\mathbf{H}_{1}\left(
t\right) }{C\left( t\right) },
\end{equation*}
where $\mathbf{H}_{1}\left( y\right) :=\mathbf{H}_{1}\left( \infty ,y\right)
,$ $\Lambda \left( y\right)$ can be estimated by
\begin{equation*}
\Lambda _{n}\left( y\right) =\int\limits_{a_{H}}^{y}\frac{d\mathbf{H}%
_{1,n}\left( t\right) }{C_{n}\left( t\right) },
\end{equation*}
where $\mathbf{H}_{1,n}\left( t\right) =\frac{1}{n}\sum\limits_{i=1}^{n}%
\mathbb{I}_{\left\{ Z_{i}\leq t,\text{ }\delta _{i}=1\right\} }.$ Obviously
\begin{equation*}
\Lambda _{n}\left( y\right) =\frac{1}{n}\sum\limits_{i=1}^{n}\frac{\mathbb{I}%
_{\left\{ Z_{i}\leq y,\text{ }\delta _{i}=1\right\} }}{C_{n}\left(
Z_{i}\right) }.
\end{equation*}
Following Lemma 2.1 of Chen and Day (2003) for $y\in \left[ a_{H},b\right] $
we have
\begin{equation}\label{6.3}
F_{n}\left( y\right) -F\left( y\right) =\left( 1-F\left( y\right) \right)
\left( \Lambda _{n}\left( y\right) -\Lambda \left( y\right) \right) +O\left(
\frac{\log \log n}{n}\right) \text{ a.s},
\end{equation}
and using integration by part we get
\begin{eqnarray*}
\Lambda _{n}\left( y\right) -\Lambda \left( y\right) &=&\left\{ \frac{%
\mathbf{H}_{1,n}\left( y\right) -\mathbf{H}_{1}\left( y\right) }{C\left(
y\right) }-\frac{\mathbf{H}_{1,n}\left( a_{H}\right) -\mathbf{H}_{1}\left(
a_{H}\right) }{C\left( a_{H}\right) }\right\} \\
&&+\int\limits_{a_{H}}^{y}\frac{C\left( t\right) -C_{n}\left( t\right) }{%
C_{n}\left( t\right) C\left( t\right) }d\mathbf{H}_{1,n}\left( t\right)
+\int\limits_{a_{H}}^{y}\frac{\mathbf{H}_{1,n}\left( t\right) -\mathbf{H}%
_{1}\left( t\right) }{C^{2}\left( t\right) }dC\left( t\right) \\
&=:&I_{1}\left( y\right) +I_{2}\left( y\right) +I_{3}\left( y\right) .
\end{eqnarray*}
From (35) in Liang and Alvarez (2011) we have
\begin{equation*}
\underset{y\in \left[ a_{H},b\right] }{\sup }\left\vert \mathbf{H}%
_{1,n}\left( y\right) -\mathbf{H}_{1}\left( y\right) \right\vert
=O_{\bf{P}}\left( n^{-1/2}\right) ,
\end{equation*}
which gives $\underset{y\in \left[ a_{H},b\right] }{\sup }\left\vert
I_{1}\left( y\right) +I_{3}\left( y\right) \right\vert =O_{\bf{P}}\left(
n^{-1/2}\right) $. Also $\ $by (i) and the fact that\\
$\int\limits_{a_{H}}^{y}\frac{d\mathbf{H}_{1,n}\left( t\right) }{%
C_{n}\left( t\right) C\left( t\right) }\rightarrow \int\limits_{a_{H}}^{y}%
\frac{d\mathbf{H}_{1}\left( t\right) }{C^{2}\left( t\right) }<\infty ,$ we
get $\underset{y\in \left[ a_{H},b\right] }{\sup }\left\vert I_{2}\left(
y\right) \right\vert =O_{\bf{P}}\left( n^{-1/2}\right) .$\\
Then%
\begin{equation*}
\underset{y\in \left[ a_{H},b\right] }{\sup }\left\vert \Lambda _{n}\left(
y\right) -\Lambda \left( y\right) \right\vert =O_{\bf{P}}\left( n^{-1/2}\right) ,
\end{equation*}
and this, together with (\ref{6.3}), give the result (ii).\\
(iii) We can easily verify that
\begin{equation*}
    \overline{H}_{n}(t)=\overline{F}_{n}(t)\overline{G}_{n}(t).
\end{equation*}
Hence observe that
\begin{equation}\label{GGn}
    H_{n}(t)-H(t)=\overline{G}(t)(F_{n}(t)-F(t))+\overline{F}(t)(G_{n}(t)-G(t))-(G_{n}(t)-G(t))(F_{n}(t)-F(t))
\end{equation}
Therefore from Lemma 5.5 of Liang (2011) and (ii) combining with (\ref{GGn}), we get the result of (iii).\\
It remains to prove (iv). Observe that
\begin{eqnarray}\label{6.5}
\mu _{n}-\mu  &=&\frac{1}{C_{n}\left( y\right) C\left( y\right) }\left\{
C\left( y\right) \overline{H}_{n}\left( y\right) \left( L_{n}\left( y\right)
-L\left( y\right) \right) -L\left( y\right) \overline{H}\left( y\right)
\left( C_{n}\left( y\right) -C\left( y\right) \right) \right.   \notag \\
&&\left. -L\left( y\right) C\left( y\right) \left( H_{n}\left( y\right)
-H\left( y\right) \right) \right\}.
\end{eqnarray}
Hence from (i) and Lemma 5.5 of Liang (2011) together with (\ref{6.5}), we get the result of (iv)
and the proof of Lemma \ref{LemmaA5} is finished.
$\hfill\Box$

\section*{References}
\small
\renewcommand{\baselinestretch}{1.3} \normalsize
\begin{description}
\item Attaoui, S., Laksaci, A. and  Ould Sa\"{\i}d, E. (2015). Asymptotic Results for an M-Estimator of the Regression Function for Quasi-Associated Processes, \textit{Functional Statistics and Applications}, Selected Papers from MICPS-2013, 3-28. Springer.

\item Benseradj, H. and Guessoum, Z. (2020). Strong uniform consistency rate of an M-estimator of regression function for incomplete data under
$\alpha$-mixing condition. \textit{Communications in Statistics - Theory and Methods}, DOI: 10.1080/03610926.2020.1764037.
\item Boente, G. and  Fraiman, R. (1990). Asymptotic distribution of robust estimators for nonparametric models from mixing processes. \textit{Ann Statist.} 18: 891-906.
\item Bradley, R.C. (2007). Introduction to Strong Mixing Conditions, Vols. 1, 2, and 3. \textit{Kendrick Press}, Heber City (Utah).

 \item Cai, Z. (1998). Asymptotic properties of Kaplan-Meier estimator for censored dependent data. \textit{Statist Probab Lett},  37: 381-389.
 \item Carbonez, A., Gy\"{o}rfi, L., and van der Meulen, E.C. (1995). Partitioning-estimates of a regression function under random censoring. \textit{Stat. Decis}, 13: 21-37.
\item Chen, Q. and Dai, Y. (2003). Kernel estimation of higher derivatives of density and hazard rate function for truncated and censored dependent data.\textit{ J. Acta. Math. Scientia}, 23B, 477-486.
\item Collomb, G. and H\"{a}rdle, W. (1986). Strong uniform convergence rates in robust nonparametric time series analysis and prediction: Kernel regression estimation from dependent observations. \textit{Stoch. Process. Appl.} 23: 77-89.

\item Doob, J.L. (1953). Stochastic Processes, NewYork: John Wiley \& Sons.
\item Ferraty, F. and Vieu, P. (2006). Nonparametric Functionnal Data Analysis, Theory and Practice. Springer-Verlag, New York.
\item Gijbels, I. and Wang, J.L. (1993). Strong representations of the survival function for truncated and censored data with applications, \textit{J. Multivariate Anal}, 47: 210-229.
\item Gorodetskii, V. V. (1977). On the strong mixing property for linear sequences. \textit{Theory Probab Appl}, 22: 411-413.
\item Guessoum, Z. and Ould Sa\"{\i}d, E. (2012). Central limit theorem for the kernel estimator of the regression function for censored time series, \textit{J Nonparametric stat}, 24: 379-397.
\item Guessoum, Z. and Tatachak, A. (2020) Asymptotic Results for Truncated-censored and Associated Data \textit{Sankhya B} 82: 142-164
\item H\"{a}rdle, W. (1984). Robust regression function estimation. \textit{J. Multivariate Anal.} 14: 169-180.
\item He, S., and Yang, G. (1998). Estimation of the truncation probability in the random truncation model, \textit{ Ann Stat}, 26: 1011-1027.
\item Iglesias-Pérez, C. and González-Manteiga, W. (1999). Strong representation of a generalized product-limit estimator for truncated and censored data with some applications, \textit{J. Nonparametr. Stat}. 10: 213-244.
\item La\"{\i}b, N. and Ould Sa\"{\i}d, E. (2000). A robust nonparametric estimation of the autoregression function under an ergodic hypothesis. \textit{Canad. J. Stat.} 28: 817-828.
\item Lemdani.M and Ould said, E. (2017). Nonparametric robust regression estimation for censored data, \textit{Statist Papers}, 58: 505-525.
\item Liang,  H.Y. and De U\~{n}a-\'{A}lvarez, J. (2011). Asymptotic properties of conditional quantile estimator for censored dependent observations, \textit{Ann. Inst. Statist. Math}, 63: 267-289.
\item Liang, H.Y. (2011). Asymptotic normality for regression function estimate under truncation and $\alpha$-mixing conditions, \textit{C. Statist Theory Methods}, 40: 1999-2021
\item Liang, H.Y., De U\~{n}a-\'{A}lvarez, J., Iglesias-Pérez, W. (2012) Asymptotic properties of conditional distribution estimator with truncated, censored and dependent data, \textit{Test} 21: 790-810.
\item Liang, H.Y. and Liu, A. A. (2013). Kernel estimation of conditional density with truncated, censored and dependent data. \textit{Journal of Multivariate Analysis}, 120: 40-58
\item Liang, H.Y., LI, D. and MIAO, T. (2015). Conditional Quantile Estimation with Truncated, Censored and Dependent Data \textit{Chin. Ann. Math.}
36B: 969-990
\item Masry, E. (2005) Nonparametric regression estimation for dependent functional data: asymptotic normality, \textit{Stochastic Process. Appl}. 115: 155-177.
\item Ould Sa\"{i}d, E. and Lemdani, M. (2006). Asymptotic properties of a nonparametric regression function estimator with randomly truncated data. \textit{Ann. Inst. Statist. Math.} 58: 357-378.

\item Robinson, P. M. (1984). Robust nonparametric auto-regression. \textit{Lecture Notes in Statist} 26: 247-255.
\item Rosenblatt, M. (1971). Markov processes. Structure and asymptotic behaviour. Berlin, Springer-Verlag, .
\item Silverman, B.W. (1986). Density Estimation for Statistics and Data Analysis. Chapman and Hall, London.
\item Sun, L.Q and Zhou, X. (2001). Survival function and density estimation for truncated dependent data. \textit{Statist Probab Lett}, 52: 47-57.
\item Tsai, W. Y, Jewell, N. P.and Wang, M. C. (1987). A note on the product limit estimator under right censoring and left truncation. \textit{Biometrika}, 74: 883-886.
\item Tsybakov, A. (1983). Robust estimates of a function. \textit{Problems Control Inform. Theory} 18: 39-52.

\item Volkonskii.VA and Rozanov.YA (1959). Some limit theorems for random functions. \textit{Theory Probab Appl} 4: 178-197.
\item Wang, J.F, and Liang, H.Y. (2012). Asymptotic properties for an M-estimator of the regression function with truncation and dependent data. \textit{J. Korean Stat. Soc}. 41: 351-367.
\item Withers, C. S. (1981). Conditions for linear processes to be strong mixing. \textit{Z. Wahrsch. verw. Gebiete}. 57: 477-480.
\item Zhengyan, L. and Chuanrong, L. (2010). Limit theory for mixing dependent random variables, \textit{Mathematic and its application}, V 378, China.

\end{description}
\end{document}